\documentclass[a4paper,12pt]{article} %izvzela 'twoside' fleqn,

\usepackage{color}

\usepackage{amsfonts, amsmath, amsthm, amssymb}
\usepackage[T1]{fontenc}
\usepackage[cp1250]{inputenc}
\usepackage{xcolor}
\usepackage{graphicx}
\usepackage{amssymb}
\usepackage{amsmath}
\usepackage{mathptmx}
\usepackage{helvet}
\usepackage{courier}
\usepackage{txfonts}
\usepackage{tikz} 
\usetikzlibrary{arrows}
\usepackage{type1cm}

\usepackage{verbatim}

\usepackage{graphicx}
\usepackage{epsfig,amscd,amssymb,amsxtra,amsmath,amsthm}
\usepackage{type1cm}
\usepackage[T1]{fontenc}
\usepackage{graphics}
\usepackage[mathscr]{eucal}
\usepackage[all]{xy}
\usepackage{amsmath,amscd}

\newcommand{\orbit}{\mathcal{O}^{\oplus}}

\newcommand{\tr}{\textup{tr}}

\newtheorem{theorem}{Theorem}[section]

\newtheorem{definition}[theorem]{Definition}
\newtheorem{lemma}[theorem]{Lemma}

\newtheorem{remark}[theorem]{Remark}

\newtheorem{observation}[theorem]{Observation}

\newcommand{\diam} {\mathop{\rm diam}\nolimits}

\newcommand{\Cl}  {\mathop{\rm Cl}\nolimits}

\begin{document}

\def\joinrel{\mkern-3mu}
\newcommand{\varproj}{\displaystyle \lim_{\multimapinv\joinrel-\joinrel-}}

\title{Lelek-like Fans: Endpoint-dense Continua Supporting Topologically Mixing Maps}
\author{Iztok Bani\v c, Goran Erceg, Ivan Jeli\' c, Judy Kennedy}
\date{}

\maketitle

\begin{abstract}
The Lelek fan is the only smooth fan that has a dense set of end-points. In this paper, we study non-smooth fans with this property; i.e., we construct an uncountable family of pairwise non-homeomorphic such fans. Furthermore, we prove that each of them admits a topologically mixing non-invertible mapping as well as a topologically mixing homeomorphism. 
\end{abstract}
\-

%%%%%%%%%%%%%%%%%%%%%%%%%%%%%%%%%%%%%%%%%%%%%%%%%%%%%%%%%%%%%%%%%%%%%%%%%%%%%%%%%
%%% I N T R O D U C T I O N S
\section{Introduction}
In this paper, we are interested in fans $X$ that have the property that the set of end-points of $X$ is dense in $X$. We call them Lelek-like fans. It has been proved by Bula and Oversteegen \cite{oversteegen}, and by Charatonik \cite{charatonik} (independently) that there is only one such smooth fan. It is called the Lelek fan and it was constructed first in \cite{lelek} by Lelek. It was studied intensively after its introduction and, recently, several papers were published about chaotic dynamics on the Lelek fan. For example,
\begin{enumerate}
	\item in \cite{banic2}, it was proved that on the Lelek fan $L$ there are a non-invertible map $f:L\rightarrow L$ and a homeomorphism $h:L\rightarrow L$ such that $(L,f)$ and $(L,h)$ are both transitive. 
	\item in \cite{USS}, it was proved that on the Lelek fan $L$ there are a non-invertible map $f:L\rightarrow L$ and a homeomorphism $h:L\rightarrow L$ such that $(L,f)$ and $(L,h)$ are both  topologically mixing. 
	\item in \cite{USS}, it was proved that on the Lelek fan $L$ there are a non-invertible map $f:L\rightarrow L$ and a homeomorphism $h:L\rightarrow L$ such that $(L,f)$ and $(L,h)$ are both topologically mixing as well as chaotic in the sense of Robinson but not in the sense of Devaney.
    \item in \cite{BE}, it was proved that there is a transitive function $f$ on a Cantor fan $X$ such that $\varprojlim(X, f )$ is a Lelek fan. In addition, the shift map on $\varprojlim(X, f )$ is a transitive homeomorphism.
	\item in \cite{USS,BE2,van}, it was proved that on the Lelek fan $L$ there are a non-invertible map $f:L\rightarrow L$ and a homeomorphism $h:L\rightarrow L$ such that  $(L,f)$ and $(L,h)$ are both  topologically mixing with non-zero entropy.
	\item in \cite{piotr}, it was proved that on the Lelek fan $L$ there is a homeomorphism $h:L\rightarrow L$ such that  $(L,h)$ is  topologically mixing with zero entropy.
	\item in \cite{piotr}, it was proved that on the Lelek fan $L$ there is a homeomorphism $h:L\rightarrow L$ such that  $(L,h)$ is completely scrambled and weakly topologically mixing.
	\item in \cite{piotrvan}, it was proved that for any $\alpha\in [0,\infty]$  on the Lelek fan $L$ there is a homeomorphism $h:L\rightarrow L$ such that $(L,h)$ is  topologically mixing and the entropy of $(L,h)$ equals $\alpha$.
\end{enumerate}
In this paper, we study the topologically structure of Lelek-like fans as well as the dynamical systems that are admitted by such fans. Theorem \ref{main1}, the main result of the paper, says that there is an uncountable family of pairwise non-homeomorphic Lelek-like fans each of which admits a topologically mixing non-invertible mapping as well as a topologically mixing homeomorphism.

We proceed as follows. In Section \ref{s2}, Lelek-like fans are introduced and basic results that are used later in the paper are presented. In Section \ref{s3}, we construct an uncountable family of pairwise non-homeomorphic Lelek-like fans. In Section \ref{s4}, a brief overview of Mahavier dynamical systems theory that is used in Section \ref{s5} is presented. In Section \ref{s5}, we construct an uncountable family of pairwise non-homeomorphic Lelek-like fans each of them admitting a topologically mixing non-invertible mapping as well as a topologically mixing homeomorphism. 

	\section{Non-smooth fans}\label{s2}
	This paper is about Lelek-like fans, i.e., fans that have a dense set of end-points. It turns out that most such fans are non-smooth. So, we dedicate this section to prove several properties about non-smooth fans.
	\begin{definition}
Let $(X,d)$ be a compact metric space. Then we define \emph{$2^X$} by 
$$
2^{X}=\{A\subseteq X \ | \ A \textup{ is a non-empty closed subset of } X\}.
$$
Let $\varepsilon >0$ and let $A\in 2^X$. Then we define  \emph{$N_d(\varepsilon,A)$} by $N_d(\varepsilon,A)=\bigcup_{a\in A}B(a,\varepsilon)$, where for each $x\in X$ and for each $r>0$, $B(x,r)$ denotes the open ball in $X$ with center $x$ and radius $r$.
The function \emph{$H_d:2^X\times 2^X\rightarrow \mathbb R$}, defined by
$$
H_d(A,B)=\inf\{\varepsilon>0 \ | \ A\subseteq N_d(\varepsilon,B), B\subseteq N_d(\varepsilon,A)\}
$$
for all $A,B\in 2^X$, is called \emph{the Hausdorff metric on $2^X$}. The pair $(2^X,H_d)$ is called \emph{the hyperspace of the space $(X,d)$}. 
\end{definition}
\begin{observation}
	Let $(X,d)$ be a compact metric space. The Hausdorff metric $H_d$ on $2^X$ is in fact a metric on $2^X$, 
\end{observation}

Let $(X,d)$ be a compact metric space, let $A$ be a non-empty closed subset of $X$,  and let $(A_n)$ be a sequence of non-empty closed subsets of $X$. When we write $\displaystyle A=\lim_{n\to \infty}A_n$, we mean $\displaystyle A=\lim_{n\to \infty}A_n$ in $(2^X,H_d)$. 
	\begin{definition}
 \emph{A continuum} is a non-empty compact connected metric space.  \emph{A subcontinuum} is a subspace of a continuum, which is itself a continuum.
 \end{definition}
\begin{definition}
Let $X$ be a continuum. 
\begin{enumerate}
\item The continuum $X$ is \emph{unicoherent} if for any subcontinua $A$ and $B$ of $X$ such that $X=A\cup B$,  the compactum $A\cap B$ is connected. 
\item The continuum $X$ is \emph{hereditarily unicoherent } provided that each of its subcontinua is unicoherent.
\item The continuum $X$ is a \emph{dendroid} if it is an arcwise connected, hereditarily unicoherent continuum.
\item If $X$ is homeomorphic to $[0,1]$, then $X$ is \emph{an arc}.   
\item Let $X$ be an arc. A point $x\in X$ is called \emph{an end-point of $X$} if  there is a homeomorphism $\varphi:[0,1]\rightarrow X$ such that $\varphi(0)=x$.
\item Let $X$ be a dendroid.  A point $x\in X$ is called an \emph{end-point of $X$} if for  every arc $A$ in $X$ that contains $x$, $x$ is an end-point of $A$.  The set of all end-points of $X$ is denoted by $E(X)$. 
\item The continuum $X$ is \emph{a simple triode} if it is homeomorphic to $([-1,1]\times \{0\})\cup (\{0\}\times [0,1])$.
\item Let $X$ be a simple triode. A point $x\in X$ is called \emph{the top-point} or, briefly, the \emph{top of $X$} if  there is a homeomorphism $\varphi:([-1,1]\times \{0\})\cup (\{0\}\times [0,1])\rightarrow X$ such that $\varphi(0,0)=x$.
\item Let $X$ be a dendroid.  A point $x\in X$ is called \emph{a ramification-point of $X$}, if there is a simple triod $T$ in $X$ with top $x$.  The set of all ramification-points of $X$ is denoted by $R(X)$. 
\item The continuum $X$ is \emph{a  fan} if it is a dendroid with at most one ramification point $v$, which is called the top of the fan $X$ (if it exists).
\item Let $X$ be a fan.   For all points $x$ and $y$ in $X$, we define  \emph{$[x,y]$} to be the arc in $X$ with end-points $x$ and $y$, if $x\neq y$. If $x=y$, then we define $[x,y]=\{x\}$.
\item Let $X$ be a fan with top $v$. We say that that the fan $X$ is \emph{smooth} if for any $x\in X$ and for any sequence $(x_n)$ of points in $X$,
$$
\lim_{n\to \infty}x_n=x \Longrightarrow \lim_{n\to \infty}[v,x_n]=[v,x].
$$ 
A fan is non-smooth if it is not smooth. 
%\item Let $X$ be a fan.  We say that $X$ is \emph{a Cantor fan} if $X$ is homeomorphic to the continuum $\bigcup_{c\in C}A_c$, where $C\subseteq [0,1]$ is the standard Cantor set and for each $c\in C$, $A_c$ is the  {convex} segment in the plane from $(0,0)$ to $(c,1)$.
\item Let $X$ be a fan.  We say that $X$ is \emph{a Lelek fan} if it is smooth and $\Cl(E(X))=X$. See Figure \ref{figure2}.
\begin{figure}[h!]
	\centering
		\includegraphics[width=20em]{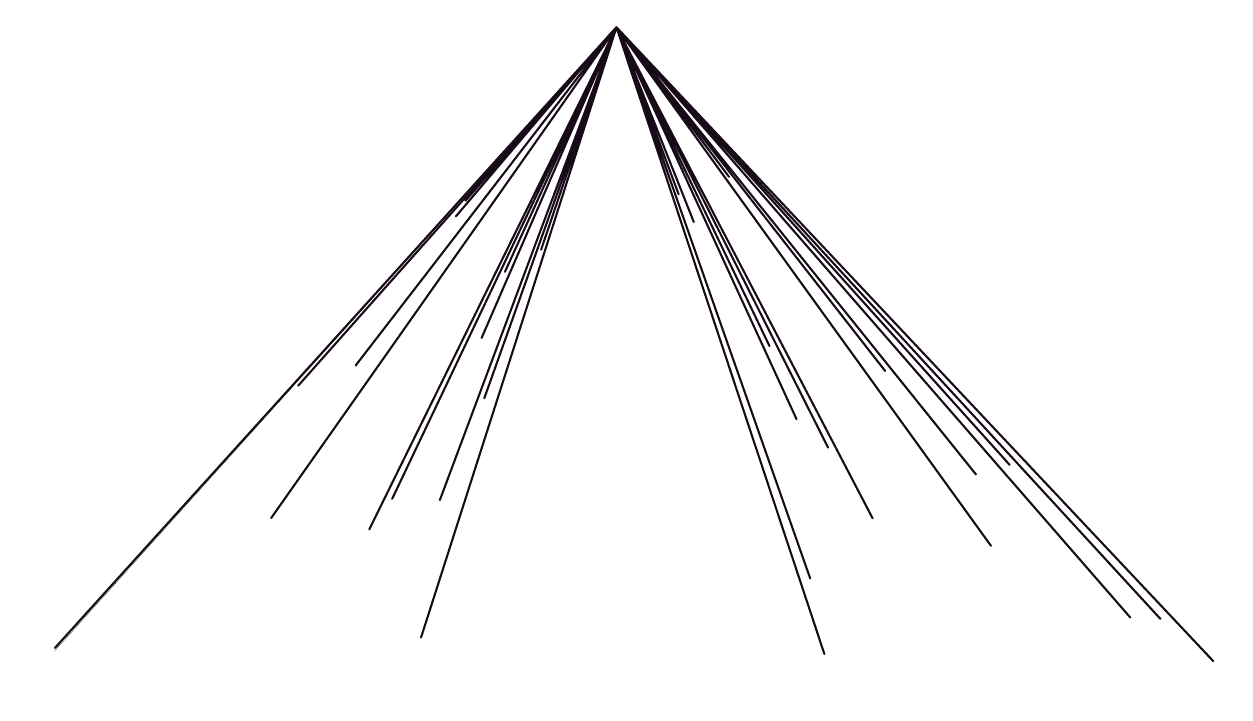}
	\caption{A Lelek fan}
	\label{figure2}
\end{figure}  
\end{enumerate}
\end{definition}
\begin{observation} \label{tatoo}
It is a well-known fact that the Lelek fan is universal for smooth fans, {i.e., every smooth fan embeds into it} (for more information see \cite{short,Jcharatonik,koch,eberhart}).   
\end{observation}
Theorem \ref{c1} was proved by Borsuk in \cite{borsuk1}.

\begin{theorem} \label{c1} Let $X$ be a fan with the top $v$. Then there is a family of arcs $\mathcal L$ in $X$ such that
\begin{enumerate}
\item\label{unk} $\displaystyle X=\bigcup_{L\in \mathcal L}L$;
\item\label{dunk} for all $L_1,L_2\in \mathcal L$,
$$
L_1\neq L_2 ~~~ \Longrightarrow ~~~  L_{1}\cap L_{2}=\{v\}.
$$
\end{enumerate}
\end{theorem}

\begin{definition}
	Let $X$ be a fan with the top $v$. Then we use $\mathcal L(X)$ to denote the family of arcs in $X$ such that
\begin{enumerate}
\item $\displaystyle X=\bigcup_{L\in \mathcal L(X)}L$;
\item for all $L_1,L_2\in \mathcal L(X)$,
$$
L_1\neq L_2 ~~~ \Longrightarrow ~~~  L_{1}\cap L_{2}=\{v\}.
$$
\end{enumerate}
The arcs in $\mathcal L(X)$ are called \emph{the legs} of the fan $X$.
\end{definition}
The following notation will be used later; it was introduced in \cite{IJGI}. 

\begin{definition}
	Let $X$ be a continuum, let $v\in X$ and let $\mathcal L$ be a collection of arcs in $X$. We say that $(v,\mathcal L)$ is \emph{a ramification pair in $X$} if 
	\begin{enumerate}
	\item $\displaystyle \bigcup_{L\in \mathcal L}L=X$, and
		\item for all $L_1,L_2\in \mathcal L$,
		$$
		L_1\neq L_2 ~~~ \Longrightarrow ~~~ L_{1}\cap L_{2}=\{v\}.
		$$
	\end{enumerate}
\end{definition}

\begin{definition}
	Let $X$ be a continuum and let $(v,\mathcal L)$ be a ramification pair in $X$. Then we say that $X$ is \emph{a pan with respect to $(v,\mathcal L)$}. We say that $X$ is \emph{a pan} if  there is a ramification pair $(v,\mathcal L)$ in $X$ such that $X$ is a pan with respect to $(v,\mathcal L)$. 
\end{definition}
\begin{definition}
	Let $X$ be a pan, let $v\in X$ and let $\mathcal L$ be a collection of arcs in $X$. We say that $X$ is a \emph{Carolyn pan with respect to $(v,\mathcal L)$} if $X$ is a pan with respect to $(v,\mathcal L)$ and for each continuum $C$ in $X$, 
	$$
	v\not \in C ~~~  \Longrightarrow  ~~~ \textup{there is } L\in \mathcal L \textup{ such that } C\subseteq L.
	$$
	We say that $X$ is \emph{a Carolyn pan} if there are a point $v\in X$ and a family $\mathcal L$ of arcs in $X$ such that $X$ is a Carolyn pan with respect to $(v,\mathcal L)$. 
%	We say that $X$ is a Carolyn phan, if $X$ is a phan as well as a Carolyn pan. If $X$ is a Carolyn phan as well as a pan with respect to $(v,\mathcal L)$, then we say that $X$ is a Carolyn phan with respect to $(v,\mathcal L)$
\end{definition}

\begin{observation}\label{lubicamoja}
	Note that each fan is a Carolyn pan.
\end{observation}
We use $\dim (X)$ to refer to the small inductive dimension \cite[Definition 1.1.1., page 3]{nagata}. It is also called the (topological) dimension of topologically spaces. Note that all fans are dendroids and all dendroids are $1$-dimensional continua \cite[(48), page 239]{jjc1}. The following is proved in \cite{IJGI}.
\begin{theorem}\label{dvojica}
	Let $X$ be a 1-dimensional pan. If $X$ is a Carolyn pan, then $X$ is a fan.
	\end{theorem}
	\begin{proof}
		See \cite[Theorem 4.28]{IJGI}.
	\end{proof}

In the remainder of this section, some new results about non-smooth fans are presented.

\begin{definition}
	Let $X$ be a fan with top $v$ and let $x\in X$. We say that $x$ is \emph{a non-smooth point in X} if there is a sequence $(x_n)$ in $X$ such that
	\begin{enumerate}
		\item for all positive integers $m$ and $n$,
		$$
		m\neq n ~~~ \Longrightarrow ~~~ x_m\neq x_n,  
		$$
		\item $\displaystyle\lim_{n\to \infty}x_n=x$, and
		\item $\displaystyle\lim_{n\to \infty}[v,x_n]\neq [v,x]$.
	\end{enumerate} 
	We use $\mathbf{NS}(X)$ to denote the set
	$$
	\mathbf{NS}(X)=\{x\in X \ | \ x \textup{ is a non-smooth point in X}\}.
	$$
\end{definition}
Lemma \ref{likija} is used in the proof of Theorem \ref{FF}.
\begin{lemma}\label{likija}
	Let $Q$ be the Hilbert cube, let $F$ be a fan in $Q$ with top $v$ and let $\mathcal A$ be a countable (countably infinite or finite) family of arcs in $Q$ such that 
	\begin{enumerate}
		\item for each $A\in \mathcal A$, $v\in E(A)$,
		\item for all $A_1,A_2\in \mathcal A$,
		$$
		A_1\neq A_2 ~~~ \Longrightarrow ~~~ A_1\cap A_2=\{v\},
		$$
		\item for each $A\in \mathcal A$, $A\cap F=\{v\}$, and
		\item if $\mathcal A=\{A_1,A_2,A_3,\ldots\}$ is countably infinite, then
		$$
		\limsup A_n\subseteq F.
		$$
	\end{enumerate} 
	Then $F\cup \left(\bigcup_{A\in \mathcal A} A\right)$ is a fan.
\end{lemma}
\begin{proof}
Let $X=F\cup \left(\bigcup_{A\in \mathcal A} A\right)$. Note that if $\mathcal A$ is finite, then $X$ is a fan. For the rest of the proof assume that $\mathcal A$ is countably infinite and let $\mathcal A=\{A_1,A_2,A_3,\ldots\}$. Note that $X$ is a union of a family of connected metric spaces all of them having a common point, $v$. Therefore, $X$ is a connected metric space. To prove that $X$ is compact, let $\mathcal O$ be a family of open sets in $Q$ such that $X\subseteq \bigcup \mathcal O$. It follows that $F\subseteq \bigcup \mathcal O$. Since $F$ is compact, there is a finite family $\{O_1,O_2,O_3,\ldots,O_m\}\subseteq \mathcal O$ such that $F\subseteq \bigcup_{i=1}^mO_i$.  It follows from $\displaystyle \limsup A_n\subseteq F$ that there is a positive integer $n_0$ such that for each positive integer $n$,
$$
n\geq n_0 ~~~ \Longrightarrow ~~~ A_n\subseteq \bigcup_{i=1}^mO_i.
$$
Since $\bigcup_{i=1}^{n_0}A_i$ is compact, there is a finite family $\{U_{1},U_2,U_3,\ldots,U_k\}\subseteq \mathcal O$ such that $\bigcup_{i=1}^{n_0}A_i\subseteq \bigcup_{i=1}^kU_i$.  Therefore,
$$
X\subseteq \left(\bigcup_{i=1}^mO_i\right)\cup \left(\bigcup_{i=1}^kU_i\right)
$$
and it follows that $X$ is compact. Thus, it is a continuum. Note that since $F$ is a fan, it follows that $F$ is either a 1-dimensional pan or $F=\{v\}$. Since (in both cases) $X$ is a union of countably many 1-dimensional continua, it follows from \cite[Theorem 4.1.9, page 257]{nagata} that it is itself a 1-dimensional continuum. Therefore, $X$ is a 1-dimensional pan. To prove that $X$ is a fan, we prove that $X$ is a Carolyn pan (this suffices by Theorem \ref{dvojica}).   
Suppose that $X$ is not a Carolyn pan. Let $C$ be a continuum in $X$ such that $v\not \in C$ and such that for each $A\in \mathcal A\cup \mathcal L(F)$, $C\not \subseteq A$. Also, let $\varepsilon >0$ be such that 
$$
\Cl(B(v,\varepsilon))\cap C=\emptyset.
$$
Note that since $F$ is a fan, it follows from Observation \ref{lubicamoja} that $F$ is a Carolyn pan. Therefore, $C\not \subseteq F$. Next, suppose that $C\subseteq \bigcup_{i=1}^{\infty}A_i$. Note that it follows that $C\cap F=\emptyset$. 
If there is a strictly increasing sequence $i_n$ of positive integers such that $A_{i_n}\cap C\neq \emptyset$, then it follows from $\limsup A_n\subseteq F$ that $C\cap F\neq \emptyset$, which is a contradiction. Therefore, there is a positive integer $n_0$ such that for each positive integer $n$,
$$
n\geq n_0 ~~~ \Longrightarrow ~~~ C\cap A_n=\emptyset.
$$
Let $F_0=\bigcup_{i=1}^{n_0}A_i$. It follows that $C\subseteq F_0$. Note that since $F_0$ is a fan, it follows from Observation \ref{lubicamoja} that $F_0$ is a Carolyn pan. Therefore, $C \subseteq A$ for some $A\in\{A_1,A_2,A_3,\ldots,A_{n_0}\}$, which is a contradiction. Therefore, $C\not\subseteq \bigcup_{i=1}^{\infty}A_i$. 
This proves that $C\cap F\neq \emptyset$ and $C\cap \bigcup_{i=1}^{\infty}A_i\neq \emptyset$. 

Let $n$ be a positive integer. We show that $C\cap A_n=\emptyset$. Suppose that $C\cap A_n\neq \emptyset$. Note that $C\cap A_n=C\cap (A_n\setminus B(v,\varepsilon))$. Then  $C\cap (A_n\setminus B(v,\varepsilon))$ and $C\cap \left(\left(\bigcup_{i\in \mathbb N\setminus\{n\}}A_i\right)\cup F\right)$ is a separation for $C$, which is a contradiction. Therefore, $C\cap A_n=\emptyset$. It follows that $C\cap \left(\bigcup_{i=1}^{\infty}A_i\right)=\emptyset$, which is a contradiction. This proves that $X$ is a Carolyn pan. 
\end{proof}
The following observation (and the notation from Definitions \ref{kulko1} and \ref{kulko2}) is used in the proof of Theorem \ref{FF}.
\begin{observation}
	For each positive integer $n$, let $a_n=n^2-n+1$. Then for each positive integer $n$,
	$$
	a_{n+1}-a_n=2n.
	$$
\end{observation}
\begin{definition}\label{kulko1}
	Let $X$ be a set and let $x\in X$. For each positive integer $n$, we use $x^n$ to denote the point $(\underbrace{x,x,x,\ldots,x}_n)\in X^n$. 
\end{definition}
\begin{definition}\label{kulko2}
	Let $Q$ be the Hilbert cube and let $\mathbf x,\mathbf y\in Q$. Then we use $I[\mathbf x,\mathbf y]$ to denote the arc
	$$
	I[\mathbf x,\mathbf y]=\{(1-t)\cdot \mathbf x+t\cdot \mathbf y \ | \ t\in [0,1]\}.
	$$ 
\end{definition}
\begin{theorem}\label{FF}
	Let $F$ be a non-degenerate fan. Then there is a fan $X$ such that $F\subseteq X$ and $$\mathbf{NS}(X)=F.$$ 
\end{theorem}
\begin{proof}
	First, suppose that $F$ has infinitely many legs and suppose that $v$ is the top of $F$. Note that every fan is a 1-dimensional continuum \cite[(48), page 239]{jjc1} and that each 1-dimensional continuum may be embedded into $\mathbb R^3$ \cite[Theorem 1.11.4, page 120]{nagata}. So, we may assume that our fan $F$ is a subspace of $[0,1]^3$ such that $v\neq (0,0,0)$. Let
	$$
	\mathbf F=F\times\{0\}\times \Big\{\big((\underbrace{0^3,0^3}_{2\cdot 1}),(\underbrace{0^3,0^3,0^3,0^3}_{2\cdot 2}), \ldots,(\underbrace{0^3,0^3,0^3,0^3,0^3,0^3,\ldots,0^3,0^3}_{2\cdot n}),\ldots\big)\Big\}.
	$$
	Note that $\mathbf F$ is homeomorphic to $F$ and that 
	$$
	\mathbf F\subseteq [0,1]^3\times\{0\}\times \Big\{\big((\underbrace{0^3,0^3}_{2\cdot 1}),(\underbrace{0^3,0^3,0^3,0^3}_{2\cdot 2}),  \ldots,(\underbrace{0^3,0^3,0^3,0^3,0^3,0^3,\ldots,0^3,0^3}_{2\cdot n}),\ldots\big)\Big\},
	$$
	which is a subspace of 
	 $$
	 [0,1]^3\times[0,1]\times \Big(\prod_{k=1}^{2\cdot 1}[0,1]^3\times \prod_{k=1}^{2\cdot 2}[0,1]^3\times \ldots \times \prod_{k=1}^{2\cdot n}[0,1]^3\times \ldots \Big)
	 $$
	and this is a topologically product of countably many closed unit intervals and, therefore, it is a copy of the Hilbert cube $\prod_{n=1}^{\infty}[0,1]$. We use $\mathbf Q$ to denote the Hilbert cube 
	$$
	\mathbf Q=[0,1]^3\times[0,1]\times \Big(\prod_{k=1}^{2\cdot 1}[0,1]^3\times \prod_{k=1}^{2\cdot 2}[0,1]^3\times \ldots \times \prod_{k=1}^{2\cdot n}[0,1]^3\times \ldots \Big).
		$$
	 Let $\{c_k \ | \ k \textup{ is a positive integer}\}$ be a countable dense subset of $F\setminus \{v\}$. For each positive integer $k$, let $L_k\in \mathcal L(F)$ be such that $c_k\in L_k$. Then $\Cl(\bigcup_{k=1}^{\infty}L_k)=F$. Also, for each positive integer $k$, let $a_k=k^2-k+1$, let $e_k=E(L_k)\setminus \{v\}$ and let 
	  $$
	  \mathbf e_k=\Big((e_k,0,\big((\underbrace{0^3,0^3}_{2\cdot 1}),(\underbrace{0^3,0^3,0^3,0^3}_{2\cdot 2}), \ldots,(\underbrace{0^3,0^3,0^3,0^3,0^3,0^3,\ldots,0^3,0^3}_{2\cdot k}),\ldots\big)\Big).
	  $$
	  Also, let
	  $$
	 \mathbf v =\Big((v,0,\big((\underbrace{0^3,0^3}_{2\cdot 1}),(\underbrace{0^3,0^3,0^3,0^3}_{2\cdot 2}), \ldots,(\underbrace{0^3,0^3,0^3,0^3,0^3,0^3,\ldots,0^3,0^3}_{2\cdot k}),\ldots\big)\Big)
	 $$
	 and for each positive integer $k$, let
	 $$
	 \mathbf L_k=\left\{\Big((x,0,\big((\underbrace{0^3,0^3}_{2\cdot 1}),(\underbrace{0^3,0^3,0^3,0^3}_{2\cdot 2}), \ldots,(\underbrace{0^3,0^3,0^3,0^3,0^3,0^3,\ldots,0^3,0^3}_{2\cdot k}),\ldots\big)\Big) \ \Big| \ x\in L_k\right\}
	 $$
	 Note that for each positive integer $k$, $\mathbf L_k\in \mathcal L(\mathbf F)$ and that $\mathbf v$ and $\mathbf e_k$ are the end-points of the leg $\mathbf L_k$.  Next, for each positive integer $n$, let $a_n=n^2-n+1$ and for each $k\in \{1,2,3,\ldots,n\}$, let 
	  $$
	  \mathbf e_k^{n,1}= \Big((e_k,\frac{1}{a_n+2k-2},\big((\underbrace{0^3,0^3}_{2\cdot 1}),(\underbrace{0^3,0^3,0^3,0^3}_{2\cdot 2}), \ldots,(\underbrace{\underbrace{0^3,0^3,0^3,\ldots,0^3}_{2k-2},e_k,\underbrace{0^3,0^3,0^3,\ldots,0^3}_{2 n-2k+1}}_{2\cdot n}),\ldots\big)\Big),
	  $$ 
	  $$
	  \mathbf e_k^{n,2}= \Big((e_k,\frac{1}{a_n+2k-1},\big((\underbrace{0^3,0^3}_{2\cdot 1}),(\underbrace{0^3,0^3,0^3,0^3}_{2\cdot 2}), \ldots,(\underbrace{\underbrace{0^3,0^3,0^3,\ldots,0^3}_{2k-1},e_k,\underbrace{0^3,0^3,0^3,\ldots,0^3}_{2 n-2k}}_{2\cdot n}),\ldots\big)\Big),
	  $$ 
	  $$
	  \mathbf f_k^{n,1}= \Big((v,\frac{1}{a_n+2k-2},\big((\underbrace{0^3,0^3}_{2\cdot 1}),(\underbrace{0^3,0^3,0^3,0^3}_{2\cdot 2}), \ldots,(\underbrace{\underbrace{0^3,0^3,0^3,\ldots,0^3}_{2k-2},v,\underbrace{0^3,0^3,0^3,\ldots,0^3}_{2 n-2k+1}}_{2\cdot n}),\ldots\big)\Big),
	  $$ 
	  and
	  $$
	  \mathbf f_k^{n,2}= \Big((v,\frac{1}{a_n+2k-1},\big((\underbrace{0^3,0^3}_{2\cdot 1}),(\underbrace{0^3,0^3,0^3,0^3}_{2\cdot 2}), \ldots,(\underbrace{\underbrace{0^3,0^3,0^3,\ldots,0^3}_{2k-1},v,\underbrace{0^3,0^3,0^3,\ldots,0^3}_{2 n-2k}}_{2\cdot n}),\ldots\big)\Big).
	  $$ 
	  Explicitly, for each positive integer $n$, we have just defined
	  $$
	  \mathbf e_1^{n,1}= \Big((e_1,\frac{1}{a_n},\big((\underbrace{0^3,0^3}_{2\cdot 1}),(\underbrace{0^3,0^3,0^3,0^3}_{2\cdot 2}), \ldots,(\underbrace{e_1,0^3,0^3,0^3,0^3,0^3,\ldots,0^3,0^3}_{2\cdot n}),\ldots\big)\Big),
	  $$ 
	  $$
	  \mathbf e_1^{n,2}= \Big((e_1,\frac{1}{a_n+1},\big((\underbrace{0^3,0^3}_{2\cdot 1}),(\underbrace{0^3,0^3,0^3,0^3}_{2\cdot 2}), \ldots,(\underbrace{0^3,e_1,0^3,0^3,0^3,0^3,\ldots,0^3,0^3}_{2\cdot n}),\ldots\big)\Big),
	  $$ 
	  $$
	  \mathbf f_1^{n,1}= \Big((v,\frac{1}{a_n},\big((\underbrace{0^3,0^3}_{2\cdot 1}),(\underbrace{0^3,0^3,0^3,0^3}_{2\cdot 2}), \ldots,(\underbrace{v,0^3,0^3,0^3,0^3,0^3,\ldots,0^3,0^3}_{2\cdot n}),\ldots\big)\Big),
	  $$ 
	  $$
	  \mathbf f_1^{n,2}= \Big((v,\frac{1}{a_n+1},\big((\underbrace{0^3,0^3}_{2\cdot 1}),(\underbrace{0^3,0^3,0^3,0^3}_{2\cdot 2}), \ldots,(\underbrace{0^3,v,0^3,0^3,0^3,0^3,\ldots,0^3,0^3}_{2\cdot n}),\ldots\big)\Big),
	  $$ 
	  $$
	  \mathbf e_2^{n,1}= \Big((e_2,\frac{1}{a_n+2},\big((\underbrace{0^3,0^3}_{2\cdot 1}),(\underbrace{0^3,0^3,0^3,0^3}_{2\cdot 2}), \ldots,(\underbrace{0^3,0^3,e_2,0^3,0^3,0^3,\ldots,0^3,0^3}_{2\cdot n}),\ldots\big)\Big),
	  $$ 
	  $$
	  \mathbf e_2^{n,2}= \Big((e_2,\frac{1}{a_n+3},\big((\underbrace{0^3,0^3}_{2\cdot 1}),(\underbrace{0^3,0^3,0^3,0^3}_{2\cdot 2}), \ldots,(\underbrace{0^3,0^3,0^3,e_2,0^3,0^3,\ldots,0^3,0^3}_{2\cdot n}),\ldots\big)\Big),
	  $$ 
	  $$
	  \mathbf f_2^{n,1}= \Big((v,\frac{1}{a_n+2},\big((\underbrace{0^3,0^3}_{2\cdot 1}),(\underbrace{0^3,0^3,0^3,0^3}_{2\cdot 2}), \ldots,(\underbrace{0^3,0^3,v,0^3,0^3,0^3,\ldots,0^3,0^3}_{2\cdot n}),\ldots\big)\Big),
	  $$ 
	  $$
	  \mathbf f_2^{n,2}= \Big((v,\frac{1}{a_n+3},\big((\underbrace{0^3,0^3}_{2\cdot 1}),(\underbrace{0^3,0^3,0^3,0^3}_{2\cdot 2}), \ldots,(\underbrace{0^3,0^3,0^3,v,0^3,0^3,\ldots,0^3,0^3}_{2\cdot n}),\ldots\big)\Big),
	  $$ 
	 $$
	 \vdots
	 $$
	  $$
	  \mathbf e_n^{n,1}= \Big((e_n,\frac{1}{a_n+2n-2},\big((\underbrace{0^3,0^3}_{2\cdot 1}),(\underbrace{0^3,0^3,0^3,0^3}_{2\cdot 2}), \ldots,(\underbrace{0^3,0^3,0^3,0^3,0^3,0^3,\ldots,e_n,0^3}_{2\cdot n}),\ldots\big)\Big),
	  $$ 
	  $$
	  \mathbf e_n^{n,2}= \Big((e_n,\frac{1}{a_n+2n-1},\big((\underbrace{0^3,0^3}_{2\cdot 1}),(\underbrace{0^3,0^3,0^3,0^3}_{2\cdot 2}), \ldots,(\underbrace{0^3,0^3,0^3,0^3,0^3,0^3,\ldots,0^3,e_n}_{2\cdot n}),\ldots\big)\Big).
	  $$ 
	  $$
	  \mathbf f_n^{n,1}= \Big((v,\frac{1}{a_n+2n-2},\big((\underbrace{0^3,0^3}_{2\cdot 1}),(\underbrace{0^3,0^3,0^3,0^3}_{2\cdot 2}), \ldots,(\underbrace{0^3,0^3,0^3,0^3,0^3,0^3,\ldots,v,0^3}_{2\cdot n}),\ldots\big)\Big),
	  $$ 
	  $$
	  \mathbf f_n^{n,2}= \Big((v,\frac{1}{a_n+2n-1},\big((\underbrace{0^3,0^3}_{2\cdot 1}),(\underbrace{0^3,0^3,0^3,0^3}_{2\cdot 2}), \ldots,(\underbrace{0^3,0^3,0^3,0^3,0^3,0^3,\ldots,0^3,v}_{2\cdot n}),\ldots\big)\Big).
	  $$ 
	  Next, for each positive integer $n$ and for each $k\in \{1,2,3,\ldots,n\}$, let 
	  $$
	  \mathbf L_k^{n,1}=\left\{\Big((x,\frac{1}{a_n+2k-2},\big((\underbrace{0^3,0^3}_{2\cdot 1}),(\underbrace{0^3,0^3,0^3,0^3}_{2\cdot 2}), \ldots,(\underbrace{\underbrace{0^3,0^3,0^3,\ldots,0^3}_{2k-2},x,\underbrace{0^3,0^3,0^3,\ldots,0^3}_{2 n-2k+1}}_{2\cdot n}),\ldots\big)\Big) \ \Big| \ x\in L_k\right\}
	$$
	and
	 $$
	  	\mathbf L_k^{n,2}=\left\{\Big((x,\frac{1}{a_n+2k-1},\big((\underbrace{0^3,0^3}_{2\cdot 1}),(\underbrace{0^3,0^3,0^3,0^3}_{2\cdot 2}), \ldots,(\underbrace{\underbrace{0^3,0^3,0^3,\ldots,0^3}_{2k-1},v,\underbrace{0^3,0^3,0^3,\ldots,0^3}_{2 n-2k}}_{2\cdot n}),\ldots\big)\Big) \ \Big| \ x\in L_k\right\}.
	$$	  
	  Explicitly, for each positive integer $n$, we have defined 
	  $$
	  	 \mathbf L_1^{n,1}=\left\{\Big((x,\frac{1}{a_n},\big((\underbrace{0^3,0^3}_{2\cdot 1}),(\underbrace{0^3,0^3,0^3,0^3}_{2\cdot 2}), \ldots,(\underbrace{x,0^3,0^3,0^3,0^3,0^3,\ldots,0^3,0^3}_{2\cdot n}),\ldots\big)\Big) \ \Big| \ x\in L_1\right\},
	$$
	 $$
	  	 \mathbf L_1^{n,2}=\left\{\Big((x,\frac{1}{a_n+1},\big((\underbrace{0^3,0^3}_{2\cdot 1}),(\underbrace{0^3,0^3,0^3,0^3}_{2\cdot 2}), \ldots,(\underbrace{0^3,x,0^3,0^3,0^3,0^3,\ldots,0^3,0^3}_{2\cdot n}),\ldots\big)\Big) \ \Big| \ x\in L_1\right\},
	$$
	 $$
	  	 \mathbf L_2^{n,1}=\left\{\Big((x,\frac{1}{a_n+2},\big((\underbrace{0^3,0^3}_{2\cdot 1}),(\underbrace{0^3,0^3,0^3,0^3}_{2\cdot 2}), \ldots,(\underbrace{0^3,0^3,x,0^3,0^3,0^3,\ldots,0^3,0^3}_{2\cdot n}),\ldots\big)\Big) \ \Big| \ x\in L_2\right\},
	  $$
	 $$
	  	 \mathbf L_2^{n,2}=\left\{\Big((x,\frac{1}{a_n+3},\big((\underbrace{0^3,0^3}_{2\cdot 1}),(\underbrace{0^3,0^3,0^3,0^3}_{2\cdot 2}), \ldots,(\underbrace{0^3,0^3,0^3,x,0^3,0^3,\ldots,0^3,0^3}_{2\cdot n}),\ldots\big)\Big) \ \Big| \ x\in L_2\right\},
	  $$
	  $$
	  \vdots
	  $$
	  $$
	  	 \mathbf L_n^{n,1}=\left\{\Big((x,\frac{1}{a_n+2n-2},\big((\underbrace{0^3,0^3}_{2\cdot 1}),(\underbrace{0^3,0^3,0^3,0^3}_{2\cdot 2}), \ldots,(\underbrace{0^3,0^3,0^3,0^3,0^3,0^3,\ldots,x,0^3}_{2\cdot n}),\ldots\big)\Big) \ \Big| \ x\in L_n\right\},
	  $$
	 $$
	  	 \mathbf L_n^{n,2}=\left\{\Big((x,\frac{1}{a_n+2n-1},\big((\underbrace{0^3,0^3}_{2\cdot 1}),(\underbrace{0^3,0^3,0^3,0^3}_{2\cdot 2}), \ldots,(\underbrace{0^3,0^3,0^3,0^3,0^3,0^3,\ldots,0^3,x}_{2\cdot n}),\ldots\big)\Big) \ \Big| \ x\in L_n\right\},
	  $$
	Note that
	\begin{enumerate}
		\item $\{ \mathbf L_k^{n,\ell} \ | \ n \textup{ is a positive integer}, k\in\{1,2,3,\ldots,n\}, \ell\in \{1,2\}\}$ is a family of mutually disjoint arcs in $\mathbf Q$,
		\item for each positive integer $n$, for each $k\in\{1,2,3,\ldots,n\}$ and for each $\ell\in \{1,2\}$, $\mathbf L_k^{n,\ell}$ is an arc in $\mathbf Q$ with end-points $\mathbf e_k^{n,\ell}$ and $\mathbf f_k^{n,\ell}$, and
		\item for each $k\in\{1,2,3,\ldots,n\}$  and for each $\ell\in \{1,2\}$, $\displaystyle \lim_{n\to \infty}\mathbf L_k^{n,\ell}=\mathbf L_k$. 
	\end{enumerate} 
	For each positive integer $n$, let 
	$$
	\mathbf A_n=I[\mathbf v,\mathbf f_1^{n,1}]\cup \mathbf L_1^{n,1}\cup I[\mathbf e_1^{n,1},\mathbf e_1^{n,2}]\cup  \mathbf L_1^{n,2}\cup I[\mathbf f_1^{n,2},\mathbf f_2^{n,1}] \cup \mathbf L_2^{n,1}\cup I[\mathbf e_2^{n,1},\mathbf e_2^{n,2}]\cup  \mathbf L_2^{n,2}\cup 
	$$
	$$
	I[\mathbf f_2^{n,2},\mathbf f_3^{n,1}] \cup  \mathbf L_3^{n,1}\cup I[\mathbf e_3^{n,1},\mathbf e_3^{n,2}]\cup  \mathbf L_3^{n,2}\cup I[\mathbf f_3^{n,2},\mathbf f_4^{n,1}] \cup \ldots \cup \mathbf L_n^{n,1}\cup I[\mathbf e_n^{n,1},\mathbf e_n^{n,2}]\cup  \mathbf L_n^{n,2}.
	$$  
	Note that
	\begin{enumerate}
		\item for each positive integer $n$, $\mathbf A_n$ is an arc in $\mathbf Q$ with end-points $\mathbf v$ and $\mathbf e_n^{n,2}$,
		\item for each positive integer $n$, $\mathbf A_n\cap \mathbf F=\{\mathbf v\}$,
		\item for all positive integers $m$ and $n$,
		$$
		m\neq n ~~~ \Longrightarrow ~~~ \mathbf A_m\cap \mathbf A_n=\{\mathbf v\},
		$$
		\item $\displaystyle \lim_{n\to \infty}\mathbf A_n=\mathbf F$.
	\end{enumerate} 
	Let 
	$$
	\mathbf X=\mathbf F\cup \bigcup_{n=1}^{\infty}\mathbf A_n.
	$$
	By Lemma \ref{likija}, $\mathbf X$ is a fan. Note that it follows from the construction of $\mathbf X$ that $\mathbf{NS}(\mathbf X)=\mathbf F$.
	
	Next, suppose that $F$ has finitely many legs: $\mathcal L(F)=\{Y_1,Y_2,Y_3,\ldots,Y_m\}$ for some positive integer $m$. Let $Y_0=Y_m$ and for each positive integer $k$, $L_k=Y_{k ~~~(\textup{mod} ~~~ m)}$. We use the construction of $\mathbf X$ from the above proof for this newly defined family $\{L_k \ | \ k \textup{ is a positive integer}\}$. We leave the details to the reader.
\end{proof}
\begin{observation}
	J. J. Charatonik and W. J. Charatonik proved in \cite{Jcharatonik1} that for each fan $X$, $X\setminus \mathbf{NS}(X)$ is dense in $X$. Therefore, for each fan $X$, $\mathbf{NS}(X)\neq X$. 
\end{observation}

\begin{definition}
	Let $X$ be a fan with top $v$ and let $A\in \mathcal L(X)$. For all $x,y\in A$ we define $x\leq_A y$ as follows:
	$$
	x\leq_A y ~~~ \Longleftrightarrow ~~~  \textup{ there is a homeomorphism } f:A\rightarrow [0,1] \textup{ such that } f(v)=0 \textup{ and } f(x)\leq f(y).  
	$$
\end{definition}
\begin{remark}
	Let $X$ be a fan with top $v$ and let $x,y\in X$. If $x$ and $y$ are elements of the same leg $A$, then we also write $x\leq y$ in $A$ instead of $x\leq_A y$.
	\end{remark}
\begin{theorem}\label{Darko}
	Let $X$ be a fan with top $v$, let $x_0\in X$ and let $(x_n)$ be a sequence in $X$ such that $\displaystyle\lim_{n\to\infty}x_n=x_0$. Then $\limsup [v,x_n]$ is a continuum.
	\end{theorem}
\begin{proof}
	This follows directly from \cite[Proposition 4.6]{IJGI}.
\end{proof}
\begin{theorem}\label{Petra}
	Let $X$ be a fan with top $v$, let $A\in \mathcal {L}(X)$, let $x_0\in A\setminus \{v\}$ and let $(x_n)$ be a sequence in $X$ such that $\displaystyle\lim_{n\to\infty}x_n=x_0$. Then 
	$$
	[v,x_0]\subseteq \limsup [v,x_n].
	$$
	\end{theorem}
\begin{proof}
		This follows directly from \cite[Proposition 4.13]{IJGI}.
\end{proof}
%\begin{theorem}
%	Let $X$ be a fan, let $A\in \mathcal {L}(X)$, let $x_0\in A\setminus \{v\}$ and let $(x_n)$ be a sequence in $X$ such that $\displaystyle\lim_{n\to\infty}x_n=x_0$. Then 
%	$$
%	\limsup [v,x_n]\subseteq A.
%	$$
%	\end{theorem}
%\begin{proof}
%	
%\end{proof}
\begin{theorem}\label{non-smooth-arc}
	Let $X$ be a fan with top $v$, and let $x_0\in \mathbf{NS}(X)$. Then there is an arc $B$ in $X$ such that
	$$
	x_0\in B\subseteq \mathbf{NS}(X).
	$$
	\end{theorem}
\begin{proof}
Let $A\in \mathcal L(X)$ be such that $x_0\in A$ and let $(x_n)$ be a sequence in $X$ such that
	\begin{enumerate}
		\item for all positive integers $m$ and $n$,
		$$
		m\neq n ~~~ \Longrightarrow ~~~ x_m\neq x_n,  
		$$
		\item $\displaystyle\lim_{n\to \infty}x_n=x_0$, and
		\item $\displaystyle\lim_{n\to \infty}[v,x_n]\neq [v,x_0]$.
	\end{enumerate} 
	Let $i_n$ be a strictly increasing sequence of positive integers such that $\displaystyle\lim_{n\to \infty}[v,x_{i_n}]$ does exist and $\displaystyle\lim_{n\to \infty}[v,x_{i_n}]\neq [v,x_0]$. We assume without loss of generality that for all positive integers $m$ and $n$,
		$$
		m\neq n ~~~ \Longrightarrow ~~~ [v,x_{i_m}]\cap [v,x_{i_n}]=\{v\}. 
		$$ 
		Note that by Theorem \ref{Petra}, $\displaystyle[v,x_0]\subseteq \lim_{n\to \infty}[v,x_{i_n}]$ and that by Theorem \ref{Darko}, $\displaystyle\lim_{n\to \infty}[v,x_{i_n}]$ is a continuum. Let $\displaystyle \hat B=\lim_{n\to \infty}[v,x_{i_n}]$ and let $z\in \hat B\setminus [v,x_0]$. 
		Note that $[z,x_0]\subseteq \hat B$.   We consider the following cases.
	\begin{enumerate}
		\item $z\in A$. Note that 
		$$
		[v,z]=[v,x_0]\cup[x_0,z]
		$$
		and let $z_0\in [x_0,z]\setminus \{x_0,z\}$. We let $B=[x_0,z_0]$.  Since $\hat B$ is the Hausdorff limit of the sequence of the arcs $[v,x_{i_n}]$, there are sequences $(z_n)$ and $(z_{0n})$ in $X$ such that
		\begin{enumerate}
		\item $\displaystyle \lim_{n\to\infty}z_n=z$ and $\displaystyle \lim_{n\to\infty}z_{0n}=z_0$,
		\item for each positive integer $n$, $z_n,z_{0n}\in [v,x_{i_n}]$
		\item for each positive integer $n$, $z_n<z_{0n}<x_{i_n}$ in $[v,x_{i_n}]$.
		\end{enumerate} 
		Note that $x_0<z_0<z$. Thus, $\displaystyle \lim_{n\to\infty}[v,z_{0n}]\neq [v,z_0]$. Then $z_0\in \mathbf{NS}(X)$ and so is $[x_0,z_0]\subseteq\mathbf{NS}(X)$.
		\item $z\not\in A$.  Note that 
		$$
		[z,v]\cup [v,x_0]=[z,x_0]
		$$
and that in this case $[v,x_0]\subseteq \mathbf{NS}(X)$ and we let $B=[v,x_0]$: let $z_0\in [v,x_0]$ and let $(z_n)$ and $(z_{0n})$ be sequences in $X$ such that
		\begin{enumerate}
		\item $\displaystyle \lim_{n\to\infty}z_n=z$ and $\displaystyle \lim_{n\to\infty}z_{0n}=z_0$,
		\item for each positive integer $n$, $z_n,z_{0n}\in [v,x_{i_n}]$
		\item for each positive integer $n$, $v<z_n<z_{0n}<x_{i_n}$ in $[v,x_{i_n}]$.
		\end{enumerate} 
Then $\displaystyle \lim_{n\to\infty}[v,z_{0n}]\neq [v,z_0]$. Thus, $[v,x_0]\subseteq\mathbf{NS}(X)$.
	\end{enumerate}  
\end{proof}
\begin{observation}
	Note that $\mathbf{NS}(X)$ need not be connected, see Figure\ref{fig333}, where such a fan $X$ is presented. 	
	\begin{figure}[h!]
	\centering
		\includegraphics[width=27.0em]{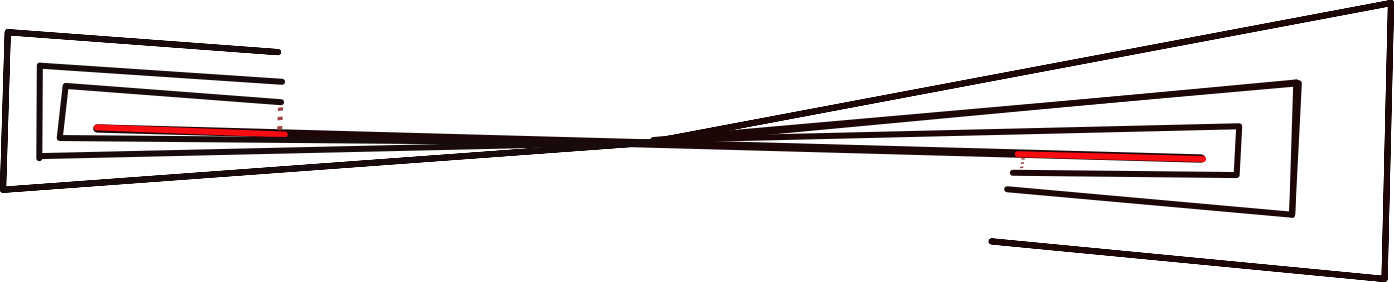}
	\caption{A fan with a non-connected non-smooth set.}
	\label{fig333}
\end{figure}

\end{observation}
\begin{definition}
	Let $X$ be a fan and let $A\in \mathcal L(X)$. We say that $A$ is a non-smooth leg in $X$, if  
	$$
	(A\cap \mathbf{NS}(X))\setminus \{v\}\neq \emptyset.
	$$
	 We use $\mathcal{NS}(X)$ to denote the set
	$$
	\mathcal{NS}(X)=\{A\in \mathcal L(X) \ | \ A \textup{ is a non-smooth leg in X}\}.
	$$ 
\end{definition}
\begin{theorem}
	Let $X$ be a fan. The following statements are equivalent.
	\begin{enumerate}
		\item\label{1} $X$ is a non-smooth fan.
		\item\label{2} $\mathbf{NS}(X)\neq \emptyset$.
		\item\label{3} $\mathcal{NS}(X)\neq \emptyset$.
	\end{enumerate}
\end{theorem}
\begin{proof}
	Note that it follows from the definition of a non-smooth fan that \ref{1} and \ref{2} are equivalent. To prove the implication from \ref{2} to \ref{3}, suppose that $\mathbf{NS}(X)\neq \emptyset$ and let $x_0\in \mathbf{NS}(X)$. Also, let $B$ be an arc in $X$ such that $x_0\in B\subseteq \mathbf{NS}(X)$. Such an arc $B$ does exist by Theorem \ref{non-smooth-arc}. It follows that $\mathbf{NS}(X)\setminus \{v\}\neq \emptyset$. Let $x\in \mathbf{NS}(X)\setminus \{v\}$ and let $A\in \mathcal L(X)$ be such that $x\in A$. It follows that $A\in \mathcal{NS}(X)$. Therefore, $\mathcal{NS}(X)\neq \emptyset$. 
   To prove the implication from \ref{3} to \ref{2}, suppose that $\mathcal{NS}(X)\neq \emptyset$ and let $A\in \mathcal{NS}(X)$. It follows that $(A\cap \mathbf{NS}(X))\setminus \{v\}\neq \emptyset$. Let $x\in (A\cap \mathbf{NS}(X))\setminus \{v\}$. Then $x\in \mathbf{NS}(X)$ and this proves that $\mathbf{NS}(X)\neq \emptyset$.
\end{proof}
\begin{definition}
		Let $X$ be a fan. We say that $X$ is \emph{a Lelek-like fan}, if $\Cl(E(X))=X$.
	\end{definition}
\begin{observation}
	Note that 
	\begin{enumerate}
		\item the Lelek fan is a Lelek-like fan.
		\item the Lelek fan is the only smooth Lelek-like fan; see \cite{oversteegen,charatonik} for details.
	\end{enumerate}
\end{observation} 
\begin{theorem}
	Let $X$ be a Lelek-like fan and let $x\in X$. The following statements are equivalent.
	\begin{enumerate}
		\item\label{11} $x\in \mathbf{NS}(X)$.
		\item\label{22} There is a sequence $(x_n)$ in $E(X)$ such that
	\begin{enumerate}
		\item for all positive integers $m$ and $n$,
		$$
		m\neq n ~~~ \Longrightarrow ~~~ x_m\neq x_n,  
		$$
		\item $\displaystyle\lim_{n\to \infty}x_n=x$, and
		\item $\displaystyle\lim_{n\to \infty}[v,x_n]\neq [v,x]$.
	\end{enumerate}
	\end{enumerate}
\end{theorem}
\begin{proof}
Note that \ref{11} follows from \ref{22}.	To prove the implication from \ref{11} to \ref{22}, suppose that $x\in \mathbf{NS}(X)$ and let $(z_n)$ be a sequence $X$ such that
	\begin{enumerate}
		\item for all positive integers $m$ and $n$,
		$$
		m\neq n ~~~ \Longrightarrow ~~~ z_m\neq z_n,  
		$$
		\item $\displaystyle\lim_{n\to \infty}z_n=x$, and
		\item $\displaystyle\lim_{n\to \infty}[v,z_n]\neq [v,x]$.
	\end{enumerate}
Since $X$ is a Lelek-like fan, it follows that for each positive integer $n$, 
$$
B\left(z_n,\frac{1}{n}\right)\cap E(X)\neq \emptyset.
$$
Let $x_1\in B(z_1,1)\cap E(X)$ and for each positive integer $n$, let 
$$
x_{n+1}\in B\left(z_{n+1},\frac{1}{n+1}\right)\cap (E(X)\setminus\{x_1,x_2,x_3,\ldots,x_n\}).
$$ 
Then $(x_n)$ is a sequence in $E(X)$ such that
	\begin{enumerate}
		\item for all positive integers $m$ and $n$,
		$$
		m\neq n ~~~ \Longrightarrow ~~~ x_m\neq x_n,  
		$$
		\item $\displaystyle\lim_{n\to \infty}x_n=x$, and
		\item $\displaystyle\lim_{n\to \infty}[v,x_n]\neq [v,x]$.
	\end{enumerate}

\end{proof}

\section{An uncountable family of Lelek-like fans}\label{s3}
In this section we prove that there is an uncountable family of pairwise non-homeomorphic Lelek-like fans. In their construction, we use quotient spaces that are defined in the following definition. 
\begin{definition}
	Let $X$ be a compact metric space and let $\sim$ be an equivalence relation on $X$. For each $x\in X$, we use $[x]$ to denote the equivalence class of the element $x$ with respect to the relation $\sim$. We also use $X/_{\sim}$ to denote the quotient space $X/_{\sim}=\{[x] \ | \ x\in X\}$. 
Also, let $q:X\rightarrow X/_{\sim}$ be the quotient map that is defined by $q(x)=[x]$ for each $x\in X$,  and let $U\subseteq X/_{\sim}$. Then 
	$$
	U \textup{ is open in } X/_{\sim} ~~~ \Longleftrightarrow ~~~ q^{-1}(U)\textup{ is open in } X.
	$$
\end{definition}

The following well-known theorem gives some basic conditions for metrizability of the quotient space $X/_{\sim}$.
	\begin{theorem}\label{tanasa}
		Let $X$ be a compact metric space and let $\sim$ be a closed equivalence relation on $X$ such that for each $x\in X$, the equivalence class $[x]$ of $x$ is a closed subset of $X$. Then the quotient space $X/_{\sim}$ is metrizable. 
	\end{theorem}
	\begin{proof}
		See \cite[Theorem 4.2.13]{engelking1}.
	\end{proof}
\begin{definition}
	We use $L$ to denote the plane dendroid $D$ (such that $E(D)$ is $1$-dimensional) which is defined in \cite[Section 9, pages 314-319]{lelek} by A.~Lelek. 
	
	We also use $v$ to denote the top of $L$. For each $A\in \mathcal L(L)$, we use $e_A$ to denote the end-point of $A$ in $A\cap E(L)$ and for each $x\in L\setminus \{v\}$, we use $A_x$ to denote the leg of $L$ that contains the point $x$.
\end{definition}
\begin{observation}
	Note that $L$ is a Lelek fan. It is the first Lelek fan ever constructed. It was shown later in \cite{charatonik}  by W.~Charatonik and  in \cite{oversteegen} by W.~D.~Bula and L.~Oversteegen  that arbitrary Lelek fans are homeomorphic. 
\end{observation}
\begin{definition}
	Let $X$ be a fan. For each $A\in \mathcal L(X)$, we choose and fix a homeomorphism $\varphi_{A,X}:A\rightarrow [0,1]$ such that $\varphi_{A,X}(v)=0$.
\end{definition}
\begin{definition}
Let $X$ be a fan, let $A\in \mathcal L(X)$, let $B$, $C$ and $D$ be arcs in $A$ and let $\varphi_{A,X}(B)=[b_1,b_2]$, $\varphi_{A,X}(C)=[c_1,c_2]$ and $\varphi_{A,X}(D)=[d_1,d_2]$. We say that the triple $(B,C,D)$ is aligned in $A$ if 
		$$
		0<b_1<b_2=c_1<c_2=d_1<d_2<1;
		$$
see Figure \ref{fig3}.
\begin{figure}[h!]
	\centering
		\includegraphics[width=17.0em]{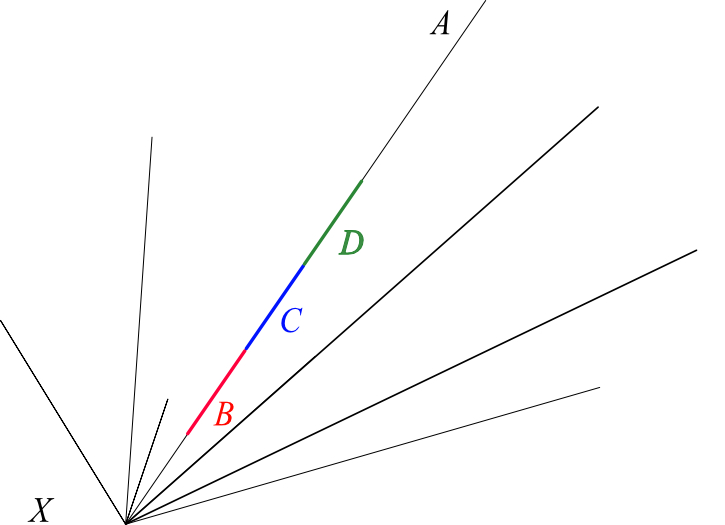}
	\caption{An aligned triple in $A$}
	\label{fig3}
\end{figure}  
\end{definition}

\begin{definition}
	Let $X$ be a fan, let $A\in \mathcal L(X)$ and let $F$ be an arc in $A$. We say that $F$ is \emph{an interruption in $A$} if for each sequence $(A_n)$ in $\mathcal L(X)$ such that $\displaystyle \lim_{n\to \infty}A_n=A$, the following holds: there is a sequence $\Big((B_n,C_n,D_n)\Big)$ of aligned triples such that for each positive integer $n$, $(B_n,C_n,D_n)$ is an aligned triple in $A_n$ and 
	$$
	\displaystyle\lim_{n\to\infty}B_n=\lim_{n\to\infty}C_n=\lim_{n\to\infty}D_n=F;
	$$
	 see Figure \ref{fig9}.
\begin{figure}[h!]
	\centering
		\includegraphics[width=24.0em]{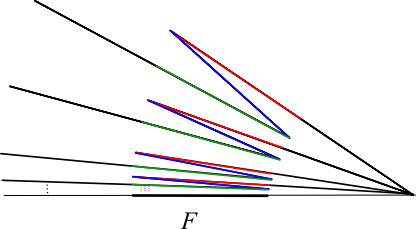}
	\caption{An interruption}
	\label{fig9}
\end{figure}    
\end{definition}
\begin{definition}
	Let $X$ be a fan, let $A\in \mathcal L(X)$, let $m$ be a positive integer, and let $F_1$, $F_2$, $F_3$, $\ldots$, $F_m$ be arcs in $A$ such that for all $i,j\in \{1,2,3,\ldots,m\}$,
	$$
	i\neq j ~~~ \Longrightarrow ~~~ F_i\cap F_j=\emptyset.
	$$ 
	We say that $\{F_1,F_2,F_3,\ldots, F_m\}$ is \emph{the set of interruptions in $A$}, if 
	\begin{enumerate}
		\item for each $i\in  \{1,2,3,\ldots,m\}$, $F_i$ is an interruption in $A$ and
		\item for each arc $F$ in $A$,
		$$
		F\not \in \{F_1,F_2,F_3,\ldots, F_m\} ~~~ \Longrightarrow ~~~ F \textup{ is not an interruption in } A.  
		$$
	\end{enumerate}
\end{definition}
\begin{definition}
	Let $X$ be a fan, let $A\in \mathcal L(X)$ and let $m$ be a positive integer. We say that \emph{the degree of $A$ is equal to $m$}, $\textup{deg}(A)=m$, if there are arcs $F_1$, $F_2$, $F_3$, $\ldots$, $F_m$ in $A$ such that the set $\{F_1,F_2,F_3,\ldots, F_m\}$ is the set of interruptions in $A$.
	
	Also, we say that \emph{the degree of $A$ is equal to $0$}, $\textup{deg}(A)=0$, if there is no set of interruptions in $A$. 
		\end{definition}
\begin{definition}
	Let $X$ be a fan, let $A\in \mathcal L(X)$, let $(B,C,D)$ be an aligned triple in $A$, and let $\varphi_{A,X}(B)=[b_1,b_2]$, $\varphi_{A,X}(C)=[c_1,c_2]$ and $\varphi_{A,X}(D)=[d_1,d_2]$. Then we define the relation $\sim_{(B,C,D)}$ on $B\cup C\cup D$ as follows: for all $x,y\in B\cup C\cup D$ we define $x \sim_{(B,C,D)} y$ if and only if one of the following hold:
	\begin{enumerate}
		\item $x=y$,
		\item $\displaystyle\varphi_{A,X}(x)=\frac{d_2-d_1}{b_2-b_1}(\varphi_{A,X}(y)-b_1)+d_1$ or $\displaystyle\varphi_{A,X}(y)=\frac{d_2-d_1}{b_2-b_1}(\varphi_{A,X}(x)-b_1)+d_1$, 
		\item $\displaystyle\varphi_{A,X}(x)=\frac{c_2-c_1}{b_1-b_2}(\varphi_{A,X}(y)-b_2)+c_1$ or $\displaystyle\varphi_{A,X}(y)=\frac{c_2-c_1}{b_1-b_2}(\varphi_{A,X}(x)-b_2)+c_1$, 
		\item $\displaystyle\varphi_{A,X}(x)=\frac{d_2-d_1}{c_1-c_2}(\varphi_{A,X}(y)-c_2)+d_1$ or $\displaystyle\varphi_{A,X}(y)=\frac{d_2-d_1}{c_1-c_2}(\varphi_{A,X}(x)-c_2)+d_1$;
	\end{enumerate}
	see Figure \ref{fig4}, where a visual presentation of the relation $\sim_{(B,C,D)}$ is presented.
\begin{figure}[h!]
	\centering
		\includegraphics[width=30.0em]{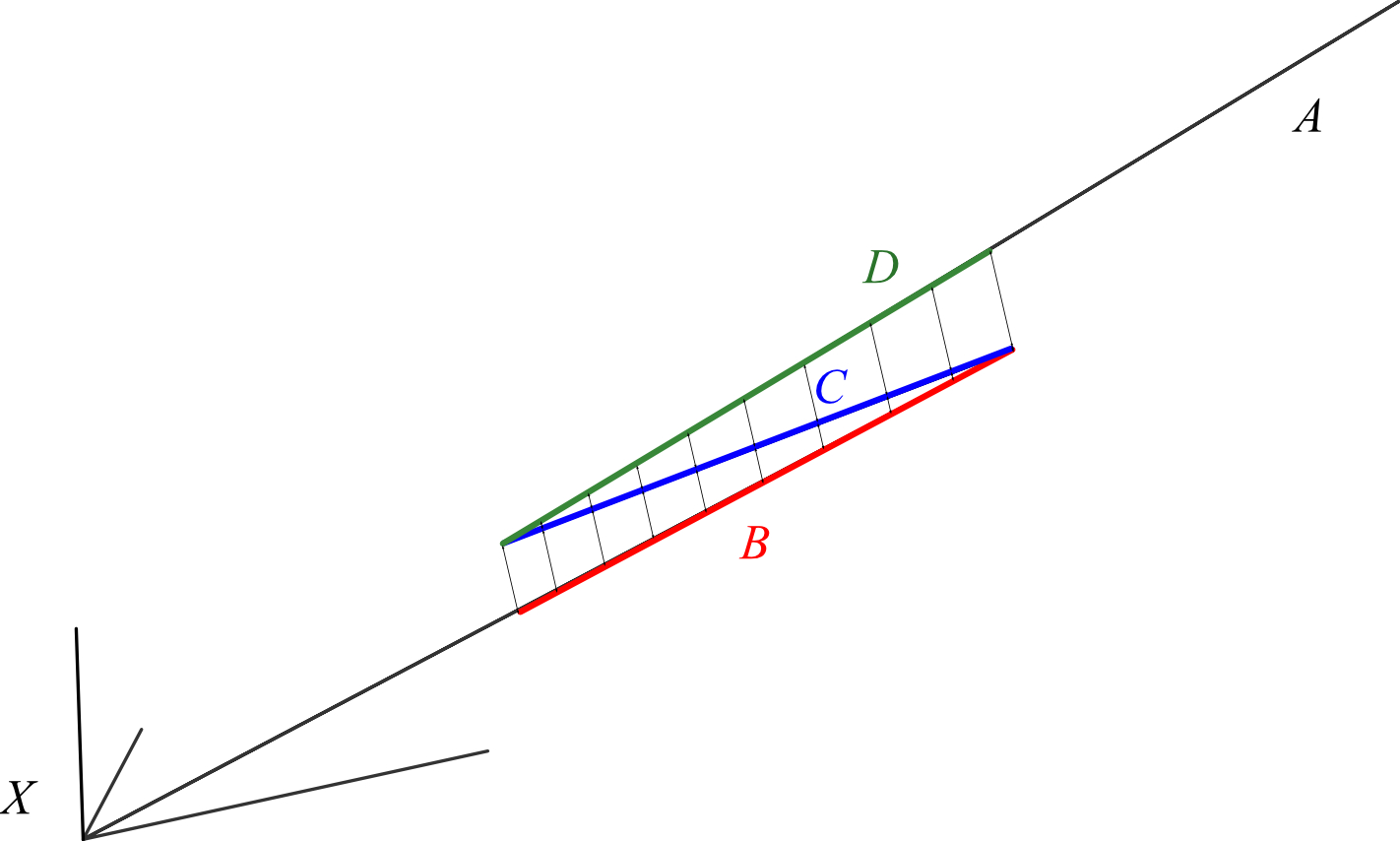}
	\caption{The relation $\sim_{(B,C,D)}$ on $B\cup C\cup D$}
	\label{fig4}
\end{figure}  
\end{definition}
\begin{observation}
Let $X$ be a fan. For each $A\in \mathcal L(X)$ and for each aligned triple $(B,C,D)$ in $A$, the relation $\sim_{(B,C,D)}$ is an equivalence relation on $B\cup C\cup D$. 
\end{observation}
\begin{definition}
	Let $X$ be a fan and let  $\sim$ be an equivalence relation on $X$. Then we use $q_{\sim}$ to denote the quotient map $q_{\sim}:L\rightarrow L/_{\sim}$, defined by 
	$$
	q_{\sim}(x)=[x]
	$$
	for each $x\in L$. 
	\end{definition}
	\begin{observation}
		Let $X$ be a fan, let $A\in \mathcal L(X)$, let $(B,C,D)$ be an aligned triple in $A$, and let $\sim$ be an equivalence relation on $X$ such that $\sim_{(B,C,D)}\subseteq \sim$. Then $$q_{\sim}(B)=q_{\sim}(C)=q_{\sim}(D).$$
	\end{observation}
	
\begin{observation}
	Let $\sim$ be an equivalence relation on $L$, let $A\in \mathcal L(L)$, and let $(B,C,D)$ be an aligned triple in $A$ such that $\sim_{(B,C,D)}\subseteq \sim$. Note that $q_{\sim}(B)$ is an interruption in $q_{\sim}(A)$; see Figure \ref{fig5}, where the topology around the interruption $q_{\sim}(B)$ in $q_{\sim}(A)$ is visualized in $L/_{\sim}$.
\begin{figure}[h!]
	\centering
		\includegraphics[width=32.0em]{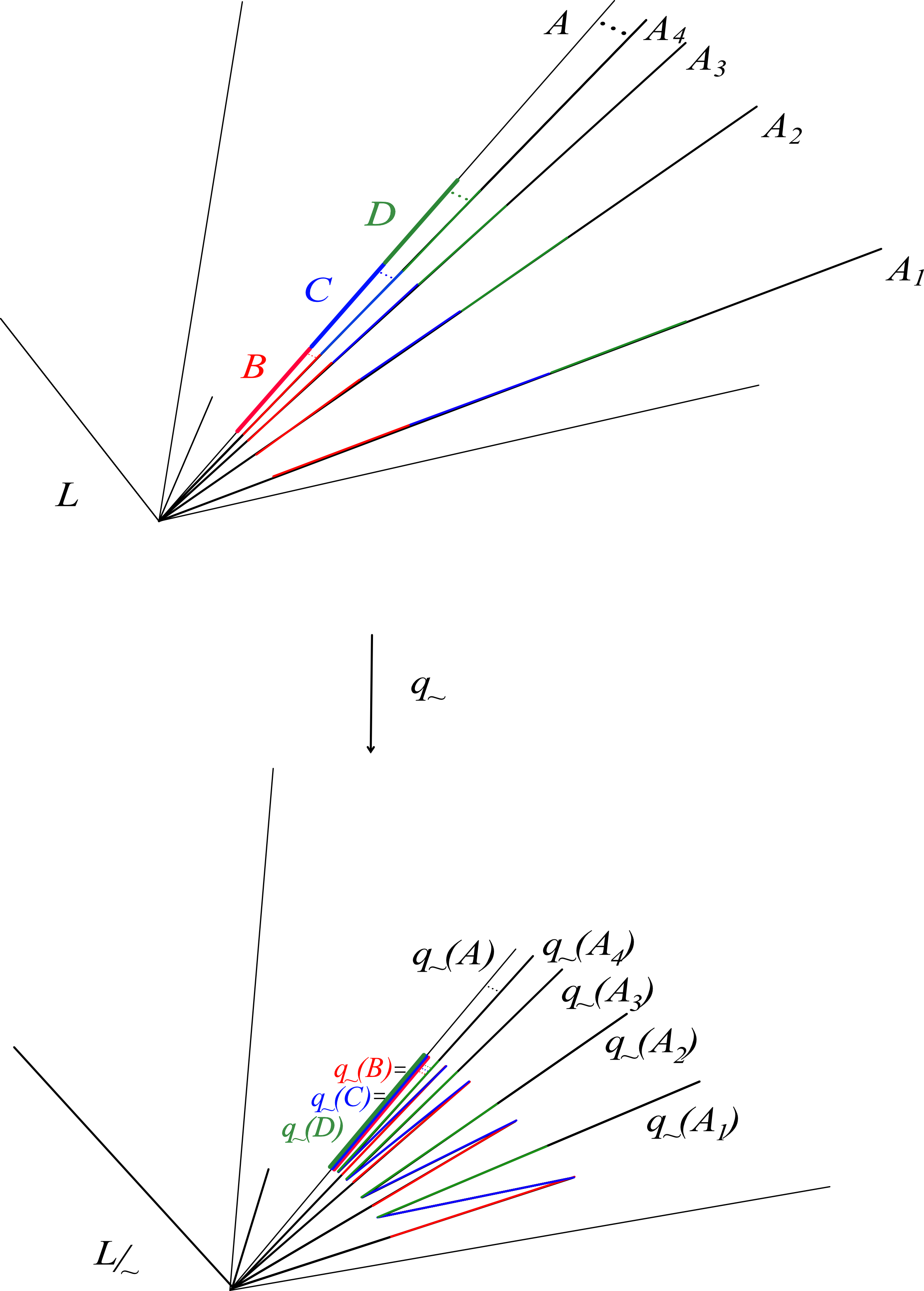}
	\caption{The interruption $q_{\sim}(B)$ in $L/_{\sim}$}
	\label{fig5}
\end{figure} 
\end{observation}
\begin{definition}
	Let $A\in \mathcal L(L)$, let $n$ be a positive integer, and for each integer  $i\in \{1,2,3,\ldots,n\}$, let $B_i$, $C_i$ and $D_i$ be arcs in $A$ such that 
	\begin{enumerate}
		\item for each $i\in \{1,2,3,\ldots,n\}$, $(B_i,C_i,D_i)$ is an aligned triple in $A$,
		\item for all $i,j\in \{1,2,3,\ldots,n\}$, 
		$$
		i\neq j ~~~  \Longrightarrow ~~~  (B_i\cup C_i\cup D_i) \cap (B_j\cup C_j\cup D_j) =\emptyset. 
		$$
	\end{enumerate} 
	Then we say that the collection $\{(B_i, C_i, D_i) \ | \ i\in \{1,2,3,\ldots,n\}\}$ is \emph{an $n$-collection of aligned triples in $A$}. 
\end{definition}
\begin{definition}
	For each positive integer $i$, let $A_i$ be a convex line segment in the plane $\mathbb R^2$ from $(0,0)$ to $\left(\frac{1}{i},\frac{1}{i^2}\right)$ and let $X$ be a continuum.  If $X$ is homeomorphic to $\bigcup_{i=1}^{\infty}A_i$, then we say that $X$ is \emph{a star}; see Figure \ref{fig6}, where a visual presentation of a star is given.
\begin{figure}[h!]
	\centering
		\includegraphics[width=25.0em]{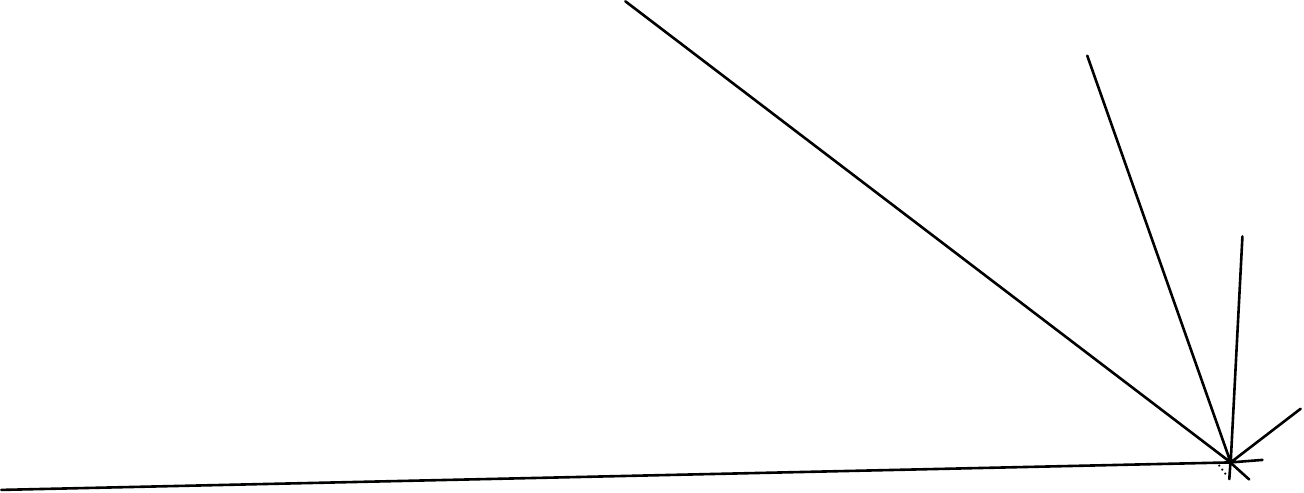}
	\caption{A star}
	\label{fig6}
\end{figure} 
\end{definition}

\begin{observation}
	Note that stars are smooth fans and since $L$ is universal for smooth fans (see Observation \ref{tatoo}), it follows that $L$ contains a star.  Furthermore, since $\Cl(E(L))=L$, it follows that there is a star $S$ in $L$ such that $\mathcal L(S)\subseteq \mathcal L(L)$.
\end{observation}
\begin{definition}
	We choose and fix a star $S$ in $L$ such that $\mathcal L(S)\subseteq \mathcal L(L)$. For each positive integer $i$, we also choose and fix $A^i\in \mathcal L(S)$ such that $S=\bigcup_{i=1}^{\infty}A^i$.
\end{definition}
\begin{observation}
Note that $\displaystyle\lim_{i\to \infty}\diam(A^i)=0$.	
\end{observation}
\begin{definition}
	We use $\mathbf I$ to denote the set
	$$
	\mathbf I=\{1,2\}\times \{3,4\}\times \{5,6\}\times  \ldots =\prod_{i=1}^{\infty}\{2i-1,2i\}.
	$$
\end{definition}
\begin{observation}\label{jur1}
	Note that $\mathbf I$ is an uncountable set.
\end{observation}
\begin{definition}
	 For each $\mathbf i=(i_1,i_2,i_3,\ldots)\in \mathbf I$ and for each positive integer $n$, we choose and fix an $i_n$-collection $\{(B_j^n, C_j^n, D_j^n) \ | \ j\in \{1,2,3,\ldots,i_n\}\}$ of aligned triples in $A^n$.
\end{definition}
\begin{definition}
	 For each $\mathbf i\in \mathbf I$ we define a relation $\sim_{\mathbf i}$ on $L$ as follows. Let $\mathbf i=(i_1,i_2,i_3,\ldots)\in \mathbf I$.  For all $x,y\in L$, we define $x\sim_{\mathbf i}y$ if and only if one of the following holds:
	 \begin{enumerate}
	 	\item $x=y$
	 	\item there is a positive integer $n$ and there is a positive integer $j\in \{1,2,3,\ldots,i_n\}$ such that $x \sim_{(B_j^n,C_j^n,D_j^n)} y$.
	 \end{enumerate}
\end{definition}
\begin{observation}
	Note that for each $\mathbf i\in \mathbf I$, the relation $\sim_{\mathbf i}$ is an equivalence relation on $L$.
\end{observation}

\begin{definition}
	Let $X$ and $Y$ be continua and let $f:X\rightarrow Y$ be a continuous function. We say that $f$ is \emph{weakly-confluent}, if for any subcontinuum $Q$ of $Y$ there is a component $C$  of $f^{-1}(Q)$ such that $f(C)=Q$.   
\end{definition}
\begin{definition}
	Let $X$ and $Y$ be continua and let $f:X\rightarrow Y$ be a continuous function. We say that $f$ is \emph{hereditarily weakly-confluent}, if for any subcontinuum $P$ of $X$ the restriction $f|_P:P\rightarrow f(P)$ is weakly confluent.   
\end{definition}
\begin{lemma}\label{confluent}
	Let $X$ and $Y$ be continua and let $f:X\rightarrow Y$ be a hereditarily weakly-confluent surjection. If $X$ is a fan, then $Y$ is a fan.
\end{lemma}
\begin{proof}
	The lemma follows from \cite[Corollary 5.23, page 149]{mackowjak}.
\end{proof}
\begin{theorem}\label{toncek}
	For each $\mathbf i\in \mathbf I$, $L/_{\sim_{\mathbf i}}$ is a fan. 
\end{theorem}
\begin{proof}
	Let $\mathbf i=(i_1,i_2,i_3,\ldots)\in \mathbf I$. By Theorem \ref{tanasa}, $L/_{\sim_{\mathbf i}}$ is metrizable, therefore, it is a continuum. Note that the quotient map $q_{\sim_{\mathbf i}}:L\rightarrow L/_{\sim_{\mathbf i}}$ is a hereditarily weakly-confluent surjection. Therefore, $L/_{\sim_{\mathbf i}}$ is a fan by Lemma \ref{confluent}.  
\end{proof}
\begin{observation}\label{endpoints}
For each $\mathbf i\in \mathbf I$, $E(L/_{\sim_{\mathbf i}})=q_{\sim_{\mathbf i}}(E(L))$.
\end{observation}
\begin{theorem}\label{jur3}
	For each $\mathbf i\in \mathbf I$, $L/_{\sim_{\mathbf i}}$ is a Lelek-like fan. 
\end{theorem}
\begin{proof}
	Let $\mathbf i=(i_1,i_2,i_3,\ldots)\in \mathbf I$. By Theorem \ref{toncek}, $L/_{\sim_{\mathbf i}}$ is a fan.  Finally, we show that $L/_{\sim_{\mathbf i}}$ is a Lelek-like fan. Let $U$ be a non-empty open set in $L/_{\sim_{\mathbf i}}$. Then $q_{\sim_{\mathbf i}}^{-1}(U)$ is a non-empty open set in $L$. Since $E(L)$ is dense in $L$, it follows that there is $e\in E(L)\cap q_{\sim_{\mathbf i}}^{-1}(U)$. Choose and fix such an end-point $e$.  It follows from Observation \ref{endpoints} that $q_{\sim_{\mathbf i}}(e)\in E(L/_{\sim_{\mathbf i}})$. Since $q_{\sim_{\mathbf i}}(e)\in U$, it follows that $E(L/_{\sim_{\mathbf i}})\cap U\neq \emptyset$. Therefore, $E(L/_{\sim_{\mathbf i}})$ is dense in $L/_{\sim_{\mathbf i}}$ and this proves that $L/_{\sim_{\mathbf i}}$ is a Lelek-like fan.
\end{proof}
\begin{observation}
	Let $\mathbf i\in \mathbf I$. Note that 
	$$
	\mathcal L(L/_{\sim_{\mathbf i}})=\{q_{\sim_{\mathbf i}}(A) \ | \ A\in \mathcal L(L)\}
	$$
	and 
	$$
	\mathcal L(L)=\{q_{\sim_{\mathbf i}}^{-1}(A) \ | \ A\in \mathcal L(L/_{\sim_{\mathbf i}})\}.
	$$
\end{observation}

	\begin{observation}
		Let $\mathbf i=(i_1,i_2,i_3,\ldots)\in \mathbf I$ and let $A\in \mathcal L(L/_{\sim_{\mathbf i}})$. Note that 
		\begin{enumerate}
		\item for each positive integer $n$, $\textup{deg}(q_{\sim_{\mathbf i}}(A^n))=i_n$ and 
		\item for each $A\in \mathcal L(L)$,  
		$$
		\textup{for each positive integer }n, ~~   A\neq A^n ~~~  \Longrightarrow  ~~~  \textup{deg}(q_{\sim_{\mathbf i}}(A))=0.
		$$
	\end{enumerate}
	\end{observation}
	\begin{lemma}\label{jure}
		Let $X$ and $Y$ be fans, let $\varphi:X\rightarrow Y$ be a homeomorphism, let $m$ be a non-negative integer, and let $A\in \mathcal L(X)$. Then the following holds
		$$
		\deg(A)=m ~~~ \Longrightarrow ~~~ \deg(\varphi(A))=m.  
		$$
	\end{lemma}
	\begin{proof}
		Let $\{F_1,F_2,F_3,\ldots,F_m\}$ be the set of interruptions in $A$. We prove that $\{\varphi(F_1),\varphi(F_2),\varphi(F_3),\ldots,\varphi(F_m)\}$ is the set of interruptions in $\varphi(A)$. Note that for all $i,j\in \{1,2,3,\ldots,m\}$,
	$$
i\neq j ~~~ \Longrightarrow ~~~ \varphi(F_i)\cap \varphi(F_j)=\emptyset.
	$$ 
 Let $(A_n)$ be any sequence of arcs in $Y$ such that $\displaystyle \lim_{n\to \infty}A_n=\varphi(A)$. Then $(\varphi^{-1}(A_n))$ is a sequence of arcs in $X$ such that $\displaystyle \lim_{n\to \infty}\varphi^{-1}(A_n)=A$. For each positive integer $n$ and for each $k\in \{1,2,3,\ldots,m\}$, let $(B_n^k,C_n^k,D_n^k)$  be an aligned triple in $\varphi^{-1}(A_n)$ such that for each $k\in \{1,2,3,\ldots,m\}$, 
		$$
		\displaystyle\lim_{n\to\infty}B_n^k=\lim_{n\to\infty}C_n^k=\lim_{n\to\infty}D_n^k=\varphi^{-1}(F_k).
		$$
		It follows that for each $k\in \{1,2,3,\ldots,m\}$, 
		$$
		\displaystyle\lim_{n\to\infty}\varphi(B_n^k)=\lim_{n\to\infty}\varphi(C_n^k)=\lim_{n\to\infty}\varphi(D_n^k)=F_k.
		$$
Finally, let $F$ be an arc in $\varphi(A)$ such that $F\not \in \{\varphi(F_1),\varphi(F_2),\varphi(F_3),\ldots,\varphi(F_m)\}$. We prove that $F$ is not an interruption in $\varphi(A)$. Since $\varphi^{-1}(F)\not \in \{F_1,F_2,F_3,\ldots,F_m\}$ and since $\{F_1,F_2,F_3,\ldots,F_m\}$ is the set of interruptions in $A$, it follows that $\varphi^{-1}(F)$ is not an interruption in $A$. Therefore, there is a sequence $(A_n)$ in $\mathcal L(X)$ such that $\displaystyle \lim_{n\to \infty}A_n=A$ and for each sequence $\Big((B_n,C_n,D_n)\Big)$ of aligned triples such that for each positive integer $n$, $(B_n,C_n,D_n)$ is an aligned triple in $A_n$, 
	$$
	\displaystyle\lim_{n\to\infty}B_n\neq F ~~~ \textup{ or } ~~~  \lim_{n\to\infty}C_n\neq F ~~~ \textup{ or } ~~~  \lim_{n\to\infty}D_n\neq F.
	$$
	Choose and fix such a sequence $(A_n)$. Suppose that $F$ is an interruption in $\varphi(A)$. Since $(\varphi(A_n))$ is a sequence in $\mathcal L(Y)$ such that $\displaystyle \lim_{n\to \infty}\varphi(A_n)=\varphi(A)$, there is a sequence $\Big((B_n,C_n,D_n)\Big)$ such that for each positive integer $n$, $(B_n,C_n,D_n)$ is an aligned triple in $\varphi(A_n)$ and 
	$$
	\displaystyle\lim_{n\to\infty}B_n=\lim_{n\to\infty}C_n=\lim_{n\to\infty}D_n=F.
	$$
	Choose and fix such a sequence $\Big((B_n,C_n,D_n)\Big)$. Then $\Big((\varphi^{-1}(B_n),\varphi^{-1}(C_n),\varphi^{-1}(D_n))\Big)$ is a sequence such that for each positive integer $n$, $(\varphi^{-1}(B_n),\varphi^{-1}(C_n),\varphi^{-1}(D_n))$ is an aligned triple in $A_n$ and 
	$$
	\displaystyle\lim_{n\to\infty}\varphi^{-1}(B_n)=\lim_{n\to\infty}\varphi^{-1}(C_n)=\lim_{n\to\infty}\varphi^{-1}(D_n)=\varphi^{-1}(F),
	$$
which is a contradiction. It follows that $F$ is not an interruption in $\varphi(A)$ and, therefore, $\{\varphi(F_1),\varphi(F_2),\varphi(F_3),\ldots,\varphi(F_m)\}$ is the set of interruptions in $\varphi(A)$.
	\end{proof}
\begin{theorem}\label{jur2}
	For all $\mathbf i,\mathbf j\in \mathbf I$, 
	$$
	\mathbf i\neq \mathbf j ~~~ \Longrightarrow ~~~  L/_{\sim_{\mathbf i}} ~~ \textup{ is not homeomorphic to } ~~ L/_{\sim_{\mathbf j}}.
	$$ 
\end{theorem}
\begin{proof}
	Let $\mathbf i,\mathbf j\in \mathbf I$ be such that $\mathbf i\neq \mathbf j$. Also, suppose that $\mathbf i=(i_1,i_2,i_3,\ldots)$ and $\mathbf j=(j_1,j_2,j_3,\ldots)$ and let $m$ be a positive integer such that $i_m\neq j_m$. Note that $\textup{deg}(q_{\sim_{\mathbf i}}(A^m))=i_m$ and that for each $A\in \mathcal L(L)$, $\textup{deg}(q_{\sim_{\mathbf j}}(A))\neq i_m$. It follows from Lemma \ref{jure} that $L/_{\sim_{\mathbf i}}$ and $L/_{\sim_{\mathbf j}}$ are not homeomorphic.         
\end{proof}
\begin{theorem}
	There is an uncountable family of pairwise non-homeomorphic Lelek-like fans.
\end{theorem}
\begin{proof}
	Let 
	$$
	\mathcal F=\{L/_{\sim_{\mathbf i}} \ | \ \mathbf i\in \mathbf I\}.
	$$ 
By Observation \ref{jur1}, 	$\mathbf I$ is uncountable and by Theorem \ref{jur3} that for each $\mathbf i\in \mathbf I$, $L/_{\sim_{\mathbf i}}$ is a Lelek-like fan. Also, it follows from Theorem  \ref{jur2} that for all $\mathbf i,\mathbf j\in \mathbf I$, 
	$$
	\mathbf i\neq \mathbf j ~~~ \Longrightarrow ~~~  L/_{\sim_{\mathbf i}} ~~ \textup{ is not homeomorphic to } ~~ L/_{\sim_{\mathbf j}}.
	$$ 
	Therefore, $\mathcal F$ is an uncountable family of pairwise non-homeomorphic Lelek-like fans.
\end{proof}
\section{Mahavier dynamical systems}\label{s4}
In Section \ref{s5}, we use Mahavier dynamical systems to construct an uncountable family of pairwise non-homeomorphic Lelek-like fans each of them admitting a topologically mixing non-invertible mapping as well as a topologically mixing homeomorphism. In this section we give an overview of the theory of Mahavier dynamical systems that is needed in Section \ref{s5}. We start with a basic dynamical system theory.

\begin{definition}
	Let $X$ be a compact metric space and let $f:X\rightarrow X$ be a continuous function. 
	We say that $(X,f)$ is \emph{a dynamical system}. 
\end{definition}
\begin{definition}
Let $(X,f)$ be a dynamical system and let $x\in X$. The sequence 
$$
\mathbf x=(x,f(x),f^2(x),f^3(x),\ldots)
$$   
is called \emph{the trajectory of $x$.} The set 
$$
\mathcal O_f^{\oplus}(x)=\{x,f(x),f^2(x),f^3(x),\ldots\}
$$   
is called \emph{the forward orbit set of $x$}.
\end{definition}
\begin{definition}
Let $(X,f)$ be a dynamical system and let $x\in X$.  If $\Cl(\orbit_f(x))=X$, then $x$ is called \emph{ a transitive point in $(X,f)$}.  We use \emph{$\tr(f)$} to denote the set
$$
\tr(f)=\{x\in X \ | \ x \textup{ is a transitive point in } (X,f)\}.
$$
\end{definition}

\begin{definition}
Let $(X,f)$ be a dynamical system.  We say that $(X,f)$ is 
\emph{transitive}  if for all non-empty open sets $U$ and $V$ in $X$,  there is a non-negative integer $n$ such that $$f^n(U)\cap V\neq \emptyset.$$ We say that the mapping $f$ is \emph{transitive}  if $(X,f)$ is transitive.
\end{definition}
The following theorem is a well-known result. See \cite{KS} for more information about transitive dynamical systems. 
\begin{theorem}\label{isolated}
Let $(X,f)$ be a dynamical system. If $X$ does not contain any isolated points, then $(X,f)$ is transitive if and only if $\tr(f)\neq \emptyset$.
\end{theorem}

%\begin{theorem}\label{andrej}
%Let $(X,f)$ be a transitive dynamical system. Then $\tr(f)$ is dense in $X$.
%\end{theorem}
\begin{definition}
Let $(X,f)$ be a dynamical system.  We say that $(X,f)$ is  \emph{  topologically mixing}  if for all non-empty open sets $U$ and $V$ in $X$,  there is a non-negative integer $n_0$ such that for each positive integer $n$,
 $$
 n\geq n_0 ~~~ \Longrightarrow ~~~  f^n(U)\cap V\neq \emptyset.
 $$
  We say that the mapping $f$ is \emph{  topologically mixing}  if $(X,f)$ is   topologically mixing. 
\end{definition}

\begin{definition}
	Let $X$ be a {non-empty} compact metric space and let ${F}\subseteq X\times X$ be a relation on $X$. If ${F}$ is closed in $X\times X$, then we say that ${F}$ is  \emph{ a closed relation on $X$}.  
\end{definition}
\begin{definition}
	Let $X$ be a {non-empty} compact metric space and let ${F}$ be a closed relation on $X$. For each positive integer $m$, we call 
	$$
	X_F^m=\Big\{(x_1,x_2,\ldots ,x_{m+1})\in \prod_{{i={1}}}^{{m{+1}}}X \ | \ \textup{for each } i\in{ \{{1,2},\ldots ,m\}}, (x_{i},x_{i+1})\in {F}\Big\}
	$$
	\emph{the $m$-th Mahavier product of ${F}$}, we call
	$$
	X_F^+=\Big\{(x_1,x_2,x_3,\ldots )\in \prod_{{i\in \mathbb N}}X \ | \ \textup{for each {positive} integer } i, (x_{i},x_{i+1})\in {F}\Big\}
	$$
	\emph{the  Mahavier product of ${F}$}, and we call
	$$
	X_F=\Big\{(\ldots,x_{-3},x_{-2},x_{-1},{x_0}{;}x_1,x_2,x_3,\ldots )\in \prod_{i\in \mathbb Z}X \ | \ \textup{for each  integer } i, (x_{i},x_{i+1})\in {F}\Big\}
	$$
	\emph{the two-sided  Mahavier product of ${F}$}.
\end{definition}

\begin{definition}\label{shit}
	Let $X$ be a {non-empty} compact metric space and let ${F}$ be a closed relation on $X$. 
	The function  $\sigma_F^{+} : {X_F^+} \rightarrow {X_F^+}$, 
	defined by 
	$$
	\sigma_F^{+} ({x_1,x_2,x_3,x_4},\ldots)=({x_2,x_3,x_4},\ldots)
	$$
	for each $({x_1,x_2,x_3,x_4},\ldots)\in {X_F^+}$, 
	is called \emph{the shift map on ${X_F^+}$}. The function  $\sigma_F : {X_F} \rightarrow {X_F}$, 
	defined by 
	$$
	\sigma_F (\ldots,x_{-3},x_{-2},x_{-1},{x_0};x_1,x_2,x_3,\ldots )=(\ldots,x_{-2},x_{-1},x_{0},{x_1};x_2,x_3,x_4,\ldots )
	$$
	for each $(\ldots,x_{-3},x_{-2},x_{-1},{x_0};x_1,x_2,x_3,\ldots )\in {X_F}$, 
	is called \emph{the shift map on ${X_F}$}.    
\end{definition}
\begin{observation}\label{juju}
	Note that $\sigma_F$ is a homeomorphism while $\sigma_F^+$ may not be a homeomorphism.
\end{observation}
\begin{definition}
	Let $X$ be a compact metric space and let $F$ be a closed relation on $X$. The dynamical system 
	\begin{enumerate}
		\item $(X_F^{+},\sigma_F^+)$ is called \emph{a Mahavier dynamical system}.
		\item $(X_F,\sigma_F)$ is called \emph{a two-sided Mahavier dynamical system}.
	\end{enumerate}
\end{definition}

We also use Theorem \ref{C} in the proof of Theorem \ref{lejgaungatamlehre}, which is one of a main results of the paper.

\begin{theorem}\label{C}
	Let $(X,G)$ be CR-dynamical system such that $p_1(G)=p_2(G)=X$ and let $F=G\cup \Delta_X$. 	Then the following hold.
	\begin{enumerate}
		\item\label{uno} If $(X_G^+,\sigma_G^+)$ is transitive, then $(X_F^+,\sigma_F^+)$ is topologically mixing.
		\item\label{due} If $(X_G,\sigma_G)$ is transitive, then $(X_F,\sigma_F)$ is topologically mixing.
	\end{enumerate}
\end{theorem}
\begin{proof}
	See \cite[Corollary 3.15]{USS}.	
\end{proof}
%
%\begin{theorem}\label{B}
%Let $X$ be a compact metric space and let $F$ be a closed relation on $X$ such that $p_1(F)=p_2(F)=X$. The following statements are equivalent. 
%\begin{enumerate}
%\item The map $\sigma_F^+$ is transitive.
%\item The homeomorphism $\sigma_F$ is transitive. 
%\end{enumerate}
%\end{theorem}
%\begin{proof}
%	The theorem follows from \cite[Theorem 4.5]{BE}.
%\end{proof}
%\begin{theorem}\label{C}
%	Let $X$ be a compact metric space and let $F$ be a closed relation on $X$. If $(X_F,\sigma_F)$ is transitive and $\Delta_X\subseteq F$, then $(X_F,\sigma_F)$ is topologically mixing.
%\end{theorem}
%\begin{proof}
%See \cite[Theorem 3.13]{USS}.
%\end{proof}	
\section{Topological mixing on Lelek-like fans}\label{s5}
In this section, we show that there is an uncountable family of pairwise non-homeomorphic Lelek-like fans each of them admitting a topologically mixing non-invertible mapping as well as a topologically mixing homeomorphism.
	We  use $I$ to denote the closed interval $I=[0,1]$. 
	\begin{definition}
	For each $(r,\rho)\in (0,\infty)\times (0,\infty)$, we define the sets \emph{$L_r$}, \emph{$L_{\rho}$} and \emph{$L_{r,\rho}$}  as follows:
${L_r}=\{(x,y)\in I\times I \ | \ y=rx\}$,  ${L_{\rho}}=\{(x,y)\in I\times I \ | \ y=\rho x\}$, and $	{L_{r,\rho}}=L_r\cup L_{\rho}$.
	%We also define the set \emph{$M_{r,\rho}$} by $M_{r,\rho}=I^+_{L_{r,\rho}}$.
\end{definition}
\begin{definition}
	Let {$(r,\rho)\in (0,\infty)\times (0,\infty)$}. We say that \emph{$r$ and $\rho$ never connect} or \emph{$(r,\rho)\in \mathcal{NC}$}  if \begin{enumerate}
		\item $r<1$, $\rho>1$ and 
		\item for all integers $k$ and $\ell$,  
		$$
		r^k = \rho^{\ell} \Longleftrightarrow k=\ell=0.
		$$
	\end{enumerate} 
\end{definition}
\begin{definition}\label{judydaljice}
	Let $(r,\rho)\in \mathcal{NC}$. We use $F_{r,\rho}$ to denote $F_{r,\rho}=L_{r,\rho}\cup\{(t,t) \ | \  t\in I\}$; see Figure \ref{fig1}.
\begin{figure}[h!]
	\centering
		\includegraphics[width=9.1em]{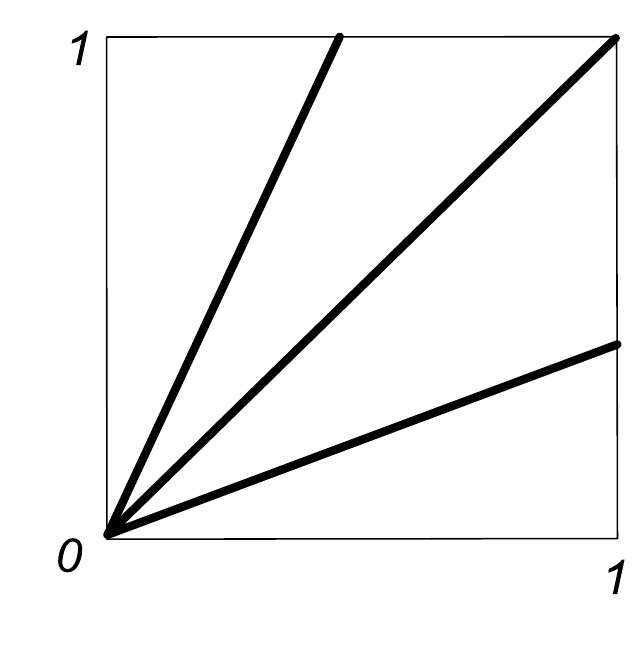}
	\caption{The relation $F_{r,\rho}$ from Definition \ref{judydaljice}}
	\label{fig1}
\end{figure}  
\end{definition}
The following two theorems are two of the main results from \cite{USS}.
\begin{theorem}
	Let $(r,\rho)\in \mathcal{NC}$. Then $I_{F_{r,\rho}}^+$ and $I_{F_{r,\rho}}$ are both  Lelek fans. 
\end{theorem}
\begin{proof}
	See the proof of \cite[Theorem 4.5, page 16]{USS}.
\end{proof}
\begin{theorem}\label{mixLelek}
	Let $(r,\rho)\in \mathcal{NC}$. The dynamical systems $(I_{F_{r,\rho}}^+,\sigma_{F_{r,\rho}}^+)$ and $(I_{F_{r,\rho}},\sigma_{F_{r,\rho}})$ are both   topologically mixing. 
\end{theorem}
\begin{proof}
	See the proof of \cite[Theorem 4.6, page 17]{USS}.
\end{proof}
\begin{definition}
	We use $\mathbb X$ to denote the set
	$$
	\mathbb X=([0,1]\cup[2,3]\cup[4,5]\cup[6,7]\cup\ldots)\cup\{\infty\}.
	$$
	\end{definition}
	We equip $\mathbb X$ with the Alexandroff one-point compactification topology $\mathcal T$; i.e., $\mathcal T$ is obtained on $\mathbb X$ from the Alexandroff one-point compactification (also known as the Alexandroff extension) of the space $[0,1]\cup[2,3]\cup[4,5]\cup[6,7]\cup\ldots$ (which is a subspace of the Euclidean line $\mathbb R$) with the point $\infty$. This topology may also be constructed as shown below by defining the metric $\textup{d}_{\mathbb X}$ on $\mathbb X$; see Definition \ref{nukja}.
	\begin{observation}\label{prej}
	For each non-negative integer $k$, let $q_k=1-\frac{1}{2^k}$ and let 
	$$
	X=([q_0,q_1]\cup [q_2,q_3]\cup [q_4,q_5]\cup [q_6,q_7]\cup \ldots )\cup \{1\}
	$$ 
	(we equip $X$ with the usual topology). Note that the compacta $\mathbb X$ and $X$ are homeomorphic. 
	\end{observation}
	\begin{definition}\label{nukja}
		Let $X$ be the compactum from Observation \ref{prej} and let $h:X\rightarrow \mathbb X$ be any homeomorphism such that for each non-negative integer $k$, $h(q_k)=k$. On the space $\mathbb X$, we always use the metric $\textup{d}_{\mathbb X}$ that is defined by 
	$
	\textup{d}_{\mathbb X}(x,y)=|h^{-1}(y)-h^{-1}(x)|
	$
	 for all $x,y\in \mathbb X$.
	\end{definition}
	\begin{observation}
	Note that the topology $\mathcal T_{\textup{d}_{\mathbb X}}$ on $\mathbb X$, that is obtained from the metric $\textup{d}_{\mathbb X}$, is exactly the one-point compactification topology $\mathcal T$ on $\mathbb X$. 	Also, note that (in this setting) for each non-negative integer $k$, 
	$	
	\diam([2k,2k+1])=\frac{1}{2^{2k+1}}.
	$
	\end{observation}
	
	\begin{definition}
		For each non-negative integer $k$, we use $I_{k+1}$ to denote  
		$
			I_{k+1}=[2k,2k+1].
			$
	\end{definition}
	\begin{observation}
		Note that for each positive integer $k$, 
		$
		\diam(I_k)=\frac{1}{2^{2k-1}}.
		$
	\end{observation}
	\begin{definition}
		We use the product metric $\textup{D}_{\mathbb X}$ on the  product $\prod_{k=-\infty}^{\infty}\mathbb X$, which is defined by 
		$$
		\textup{D}_{\mathbb X}(\mathbf x,\mathbf y)=\sup \Big\{\frac{\textup{d}_{\mathbb X}(\mathbf x(k),\mathbf y(k))}{2^{|k|}} \ \big| \ k \textup{ is an integer}\Big\}
		$$
		for all $\mathbf x,\mathbf y\in \prod_{k=-\infty}^{\infty}\mathbb X$. 
	\end{definition}
%	\begin{observation}
%		Since $\mathbb X$ is compact it follows that for all $\mathbf x,\mathbf y\in \prod_{k=-\infty}^{\infty}\mathbb X$, 
%		$$
%		\sup \Big\{\frac{\textup{d}_{\mathbb X}(\mathbf x(k),\mathbf y(k))}{2^{|k|}} \ \big| \ k \textup{ is an integer}\Big\}=\max \Big\{\frac{\textup{d}_{\mathbb X}(\mathbf x(k),\mathbf y(k))}{2^{|k|}} \ \big| \ k \textup{ is an integer}\Big\}.
%		$$
%	\end{observation}
	Next, we define the closed relation $H$ on $\mathbb X$ that will play an important role in our construction of an uncountable family of pairwise non-homeomorphic Lelek-like fans that admit topologically mixing homeomorphisms.
\begin{definition}
	Let $(r,\rho)\in \mathcal{NC}$. We use $H$ to denote the closed relation on $\mathbb X$ that is defined as follows:
	\begin{align*}
		H=&\Big\{\big(t,r\cdot t\big) \ \big| \ t\in I_1\Big\}\cup \Big\{\big(t,\rho\cdot t\big) \ \big| \ t\in I_1\Big\}\cup \\
		&\Big\{(t,t+2) \ \big| \ t\in I_1\cup I_2\cup I_3\cup I_4\cup\ldots\Big\}\cup\\
		&\Big\{(t,t-2) \ \big| \ t\in I_2\cup I_3\cup I_4\cup I_5\cup\ldots\Big\}\cup \\
		&\Big\{(t,t) \ \big| \ t\in I_1\cup I_2\cup I_3\cup I_4\cup\ldots\Big\}\cup \\
		&\Big\{(\infty,\infty)\Big\};
	\end{align*}
	see Figure \ref{uncun}. 
	\begin{figure}[h!]
	\centering
		\includegraphics[width=30em]{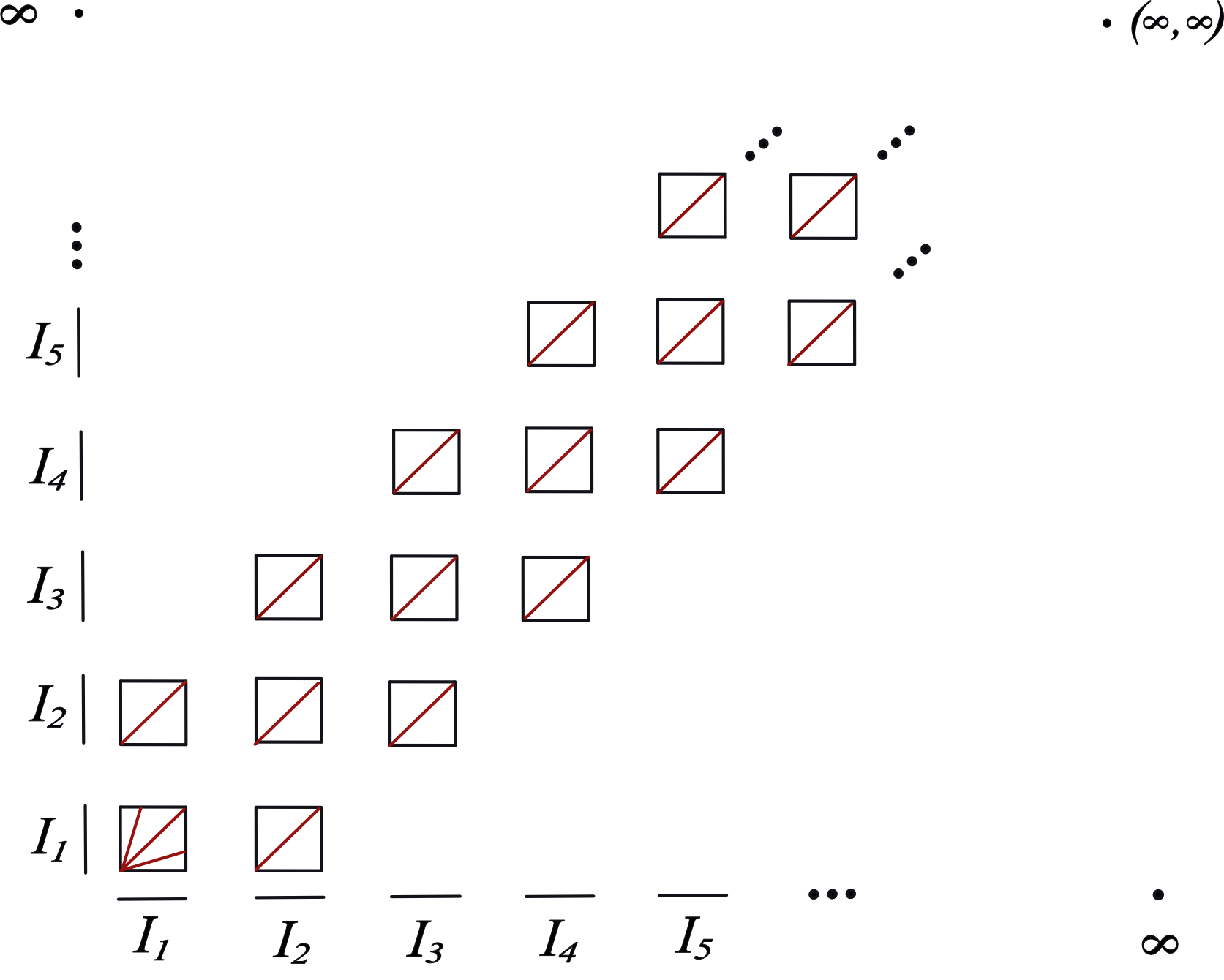}
	\caption{The relation $H$ on $\mathbb X$}
	\label{uncun}
\end{figure}
 
 \noindent We also use $\sigma_{H}^+$ to denote the shift map on the Mahavier product $\mathbb X_{H}^+$ and $\sigma_{H}$ to denote the shift map on the two-sided Mahavier product $\mathbb X_{H}$. 
\end{definition}
Next, we prove in Theorem \ref{lejgaungatamlehre} that the dynamical system $(\mathbb X_H,\sigma_H)$ is topologically mixing. In its proof, we use Lemma \ref{transitiveR}.
\begin{lemma}\label{transitiveR}
	Let $R$ be a closed relation on $I=I_1=[0,1]$ and let $H_R$ denote the closed relation on $\mathbb X$ that is defined as follows:
	\begin{align*}
		H_R=&R ~ \cup \\
		&\Big\{(t,t+2) \ \big| \ t\in I_1\cup I_2\cup I_3\cup I_4\cup\ldots\Big\}\cup \\
		&\Big\{(t,t-2) \ \big| \ t\in I_2\cup I_3\cup I_4\cup I_5\cup\ldots\Big\}\cup \\
		&\Big\{(t,t) \ \big| \ t\in I_2\cup I_3\cup I_4\cup I_5\cup\ldots\Big\}\cup \\
		& \Big\{(\infty,\infty)\Big\};
	\end{align*}
	see Figure \ref{uncun11}. 
	\begin{figure}[h!]
	\centering
		\includegraphics[width=30em]{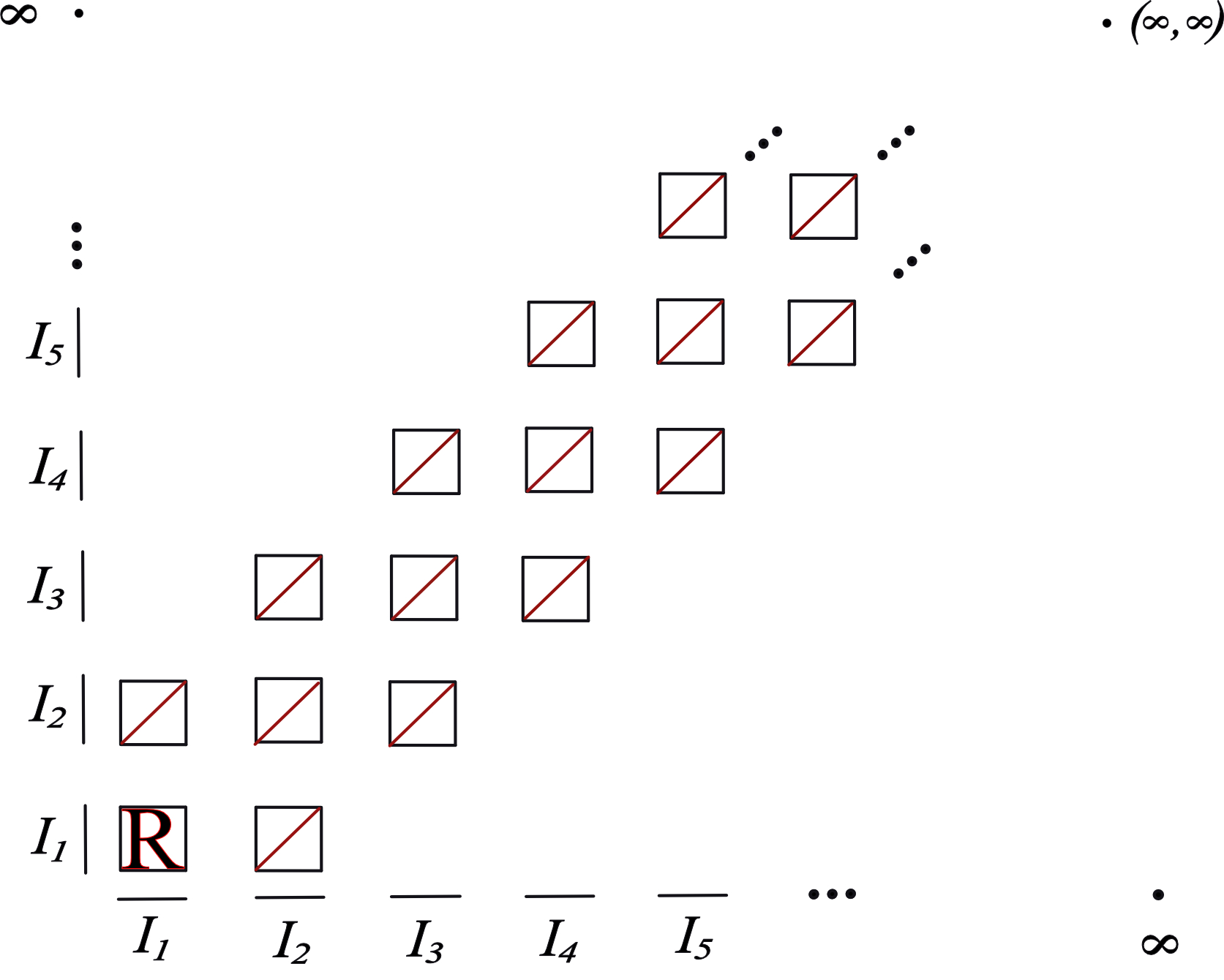}
	\caption{The relation $H_R$ on $\mathbb X$}
	\label{uncun11}
\end{figure}

\noindent Also, let $\sigma_R^+$ be the shift map on $I_R^+$, and let $\sigma_{H_R}^+$ be the shift map on $\mathbb X_{H_R}^+$. 
If $(I_R^+,\sigma_R^+)$ is transitive, then also $(\mathbb X_{H_R}^+,\sigma_{H_R}^+)$ is transitive.
\end{lemma}
\begin{proof}
	 Note that neither of the spaces $I_R^+$ or $\mathbb X_{H_R}^+$ contains an isolated point. By Theorem \ref{isolated}, $(I_R^+,\sigma_R^+)$ is transitive if and only if $\tr(\sigma_R^+)\neq \emptyset$ and  $(\mathbb X_{H_R}^+,\sigma_{H_R}^+)$ is transitive if and only if $\tr(\sigma_{H_R}^+)\neq \emptyset$. Suppose that $(I_R^+,\sigma_R^+)$ is transitive and let $(x_1,x_2,x_3,\ldots)\in \tr(\sigma_R^+)$. Then 
	 $$
	(x_1,x_1+2,x_1,x_2,x_2+2, x_2+4, x_2+2, x_2, x_3,x_3+2, x_3+4, x_3+6,x_3+4, x_3+2, x_3,x_4,\ldots) \in \tr(\sigma_{H_R}^+)
	 $$ 
	 and this proves that $\tr(\sigma_{H_R}^+)\neq \emptyset$. Therefore, $(\mathbb X_{H_R}^+,\sigma_{H_R}^+)$ is transitive.
\end{proof}

\begin{theorem}\label{lejgaungatamlehre}
	The dynamical systems $(\mathbb X_H^+,\sigma_H^+)$ and $(\mathbb X_H,\sigma_H)$ are topologically mixing.
\end{theorem}
\begin{proof}
	Since $p_1(H)=p_2(H)=\mathbb X$, it suffices to see that $(\mathbb X_H^+,\sigma_H^+)$ is transitive (by Theorem \ref{C}). Let $F=H\cap (I_1\times I_1)$ and let $I=I_1$. Then $(I_F^+,\sigma_F^+)$ is transitive by Theorem \ref{mixLelek}. It follows from Lemma \ref{transitiveR} that  $(\mathbb X_H^+,\sigma_H^+)$ is transitive. 
\end{proof}
Next, we show in Theorem \ref{lilija} that there are a topologically mixing homeomorphism $\varphi:L\rightarrow L$ and  a sequence $(A_n)$ in $\mathcal L(L)$ such that $\displaystyle \lim_{n\to \infty}\diam(A_n)=0$ and such that for each positive integer $n$ and for each $x\in A_n$,
	$$
	\varphi(x)=x.
	$$
To do that, we define  an equivalence relation $\sim$ on $\mathbb X_H$ in such a way that the quotient space $\mathbb X_H/_{\sim}$ is homeomorphic to the Lelek fan $L$. 

\begin{definition}
	We use $\sim_1$ and $\sim_2$ to denote the equivalence relations in $\mathbb X_H^+$ and $\mathbb X_H$, respectively, as follows:
	\begin{enumerate}
		\item for all $\mathbf x=(x_1,x_2,x_3,\ldots), \mathbf y=(y_1,y_2,y_3,\ldots)\in \mathbb X_H^+$, 
	$$
	\mathbf x\sim_1 \mathbf y ~~~  \Longleftrightarrow  ~~~ (\mathbf x= \mathbf y) \textup{ or (for each } k, x_k,y_k \in \{2i \ | \ i \textup{ is a non-negative integer}\}).
	$$
	\item for all $\mathbf x=(\ldots, x_{-2},x_{-1},x_0,x_1,x_2,\ldots), \mathbf y=(\ldots, y_{-2},y_{-1},y_0,y_1,y_2,\ldots)\in \mathbb X_H$, 
	$$
	\mathbf x\sim_2 \mathbf y ~~~  \Longleftrightarrow  ~~~ (\mathbf x= \mathbf y) \textup{ or (for each } k, x_k,y_k \in \{2i \ | \ i \textup{ is an integer}\}).
	$$
	\end{enumerate} 
	For each $\mathbf x\in \mathbb X_H^+$, we use $[\mathbf x]_1$ to denote the equivalence class of $\mathbf x$ with respect to relation $\sim_1$, and, for each $\mathbf x\in \mathbb X_H$, we use $[\mathbf x]_2$ to denote the equivalence class of $\mathbf x$ with respect to relation $\sim_2$. 
	\end{definition}
%In the proof of Theorem \ref{jurceksito}, we use confluent mappings that are defined in the following definition. 
%\begin{definition}
%	Let $X$ and $Y$ be any continua and let $f:X\rightarrow Y$ be a continuous mapping. We say that $f$ is \emph{confluent}, if for every subcontinuum $S$ of $Y$ and for each component $C$ of $f^{-1}(S)$, $	f(C)=S$. 
%\end{definition}
%The following is a well-known result.
%\begin{theorem}\label{char}
%	Let $X$ and $Y$ be any continua and let $f:X\rightarrow Y$ be a confluent surjection. If $X$ is a smooth fan, then also $Y$ is a smooth fan.
%	\end{theorem}
%\begin{proof}
%	See \cite[Theorem 13, page 33]{Jcharatonik}.
%\end{proof}

\begin{theorem}\label{jurceksito}
	 The quotient spaces $\mathbb X_H^+/_{\sim_1}$ and $\mathbb X_H^+/_{\sim_2}$ are both homeomorphic to the Lelek fan $L$. 
\end{theorem}
\begin{proof}
	Let $\mathbb S=\{2k \ | \ k \textup{ is a non-negative integer}\}$. First, note that $\mathbb X_H^+$ and $\mathbb X_H$  are both a $1$-dimensional compacta, both unions of smooth fans.  Each of these smooth fans is 
	\begin{enumerate}
		\item either a convex arc with one of its end-points having all coordinates in the set $\mathbb S$ 
		\item or a fan, whose legs are convex arcs, with the top having all coordinates in the set $\mathbb S$.  
	\end{enumerate} 
	Then the quotient spaces $\mathbb X_H^+/_{\sim_1}$ and $\mathbb X_H/_{\sim_2}$ are obtained from $\mathbb X_H^+$ and $\mathbb X_H$, respectively, by gluing all such points having all coordinates in the set $\mathbb S$ into one point. It follows from this that $\mathbb X_H^+/_{\sim_1}$ and $\mathbb X_H/_{\sim_2}$ are both smooth fans. In the proof of \cite[Theorem 4.5, page 16]{USS}, it is proved that for each $\mathbf t=(t_1,t_2,t_3,\ldots)\in I_{F_{r,\rho}}^+$,  
	$$
	\sup\{t_n \ | \ n \textup{ is a positive integer}\}=1 ~~~ \Longleftrightarrow ~~~ \mathbf t\in E(I_{F_{r,\rho}}^+).
	$$
	It follows from this that
	\begin{enumerate}
		\item for each $\mathbf x=(x_1,x_2,x_3,\ldots)\in \mathbb X_H^+$, it holds that if there is a positive integer $n_0$ such that $	\sup\{x_n \ | \ n\geq n_0\}=1$, 	then $[\mathbf x]_1\in E(\mathbb X_H/_{\sim_1})$. 
		\item for each $\mathbf x=(\ldots, x_{-2},x_{-1},x_0,x_1,x_2,\ldots)\in \mathbb X_H$, it holds that if there is a positive integer $n_0$ such that $	\sup\{x_n \ | \ n\geq n_0\}=1$, 	then $[\mathbf x]_2\in E(\mathbb X_H/_{\sim_2})$. 
	\end{enumerate} 
	First, we show that $\mathbb X_H^+/_{\sim_1}$ is homeomorphic to the Lelek fan $L$ by showing that the set $E(\mathbb X_H^+/_{\sim_1})$ is dense in $\mathbb X_H^+/_{\sim_1}$.  Let $\mathbf x=(x_1,x_2,x_3,\ldots)\in \mathbb X_H^+$, let $\varepsilon >0$, let $\delta>0$ be such that for all $\mathbf s,\mathbf t\in \mathbb X_H^+$, 
	$$
	d(\mathbf s,\mathbf t)<\delta ~~~ \Longrightarrow ~~~ d_{1}([\mathbf s]_1,[\mathbf t]_1)<\varepsilon
	$$
	(where $d$ is the metric on $\mathbb X_H^+$ and $d_1$ is the metric on $\mathbb X_H/_{\sim_1}$), and let $n_0$ be a positive integer such that $\frac{1}{2^{n_0}}<\delta$.  Next,
	\begin{enumerate}
		\item let $t=x_{n_0+1}$, 
		\item let $\ell$ be a non-negative integer such that $t-2\ell\in [0,1]$, 
		\item let $s=t-2\ell$, and 
		\item let $(a_n)$ be a sequence of $r$'s and $\rho$'s such that 
		$$
	\sup\{a_n\cdot a_{n-1}\cdot a_{n-2}\cdot \ldots \cdot a_3\cdot a_2\cdot a_1\cdot s \ | \ n \textup{ is a positive integer}\}=1	
		$$
		(such a sequence does exist by \cite[Theorem 9]{banic1})
	\end{enumerate} 
	and let
	$$
	\mathbf y=(x_1,x_2,x_3,\ldots,x_{n_0},t,t-2,t-4, \ldots t-2(\ell-1),s,a_1\cdot s,a_2\cdot a_1\cdot s,a_3\cdot a_2\cdot a_1\cdot s,\ldots).
	$$
	Note that $d_{1}([\mathbf x]_1,[\mathbf y]_1)<\varepsilon$ and that $[\mathbf y]_1\in  E(\mathbb X_H^+/_{\sim_1})$. 
	
	The proof that $\mathbb X_H/_{\sim_2}$ is homeomorphic to the Lelek fan $L$ is similar to the proof that $\mathbb X_H^+/_{\sim_1}$ is homeomorphic to the Lelek fan $L$; we leave the details to a reader.
	\end{proof}

We show in Theorem \ref{lijk} that $\sigma_H^{\star}$ is a homeomorphism on $\mathbb X_H/_{\sim}$ such that $(\mathbb X_H/_{\sim},\sigma_H^{\star})$ is topologically mixing. First, we recall the following basic definitions about quotient dynamical systems. 
\begin{definition}
	Let $X$ be a compact metric space, let $\sim$ be an equivalence relation on $X$,  and let $f:X\rightarrow X$ be a function such that for all $x,y\in X$,
	$$
	x\sim y  \Longleftrightarrow f(x)\sim f(y).
	$$
	Then we let $f^{\star}:X/_{\sim}\rightarrow X/_{\sim}$ be defined by   
	$
	f^{\star}([x])=[f(x)]
	$
	for any $x\in X$. 
\end{definition}
\begin{observation}\label{mutula}
	Let $(X,f)$ be a dynamical system. Note that we have defined a dynamical system as a pair of a compact metric space with a continuous function on it and that in this case, $X/_{\sim}$ is not necessarily metrizable. So, if $X/_{\sim}$ is metrizable, then also $(X/_{\sim},f^{\star})$ is a dynamical system. Note that in this case, $X/_{\sim}$ is semi-conjugate to $X$: for $\alpha:X\rightarrow X/_{\sim}$, defined by $\alpha(x)=q(x)$ for any $x\in X$, where $q$ is the quotient map obtained from $\sim$, $\alpha\circ f=f^{\star}\circ \alpha$.
\end{observation}
\begin{definition}
	Let $(X,f)$ be a dynamical system and let $\sim$ be an equivalence relation on $X$  such that for all $x,y\in X$,
	$$
	x\sim y  \Longleftrightarrow f(x)\sim f(y).
	$$
	Then we say that $(X/_{\sim},f^{\star})$ is \emph{a quotient of the dynamical system $(X,f)$} or it is \emph{the quotient of the dynamical system $(X,f)$ that is obtained from the relation $\sim$}. 
\end{definition}
We we use Theorem \ref{tupki} to prove Theorem \ref{lijk}.
\begin{theorem}\label{tupki}
	Let $X$ be a compact metric space, let $\sim$ be an equivalence relation on $X$, and  let $f:X\rightarrow X$ be a function such that for all $x,y\in X$,
$$
x\sim y  \Longleftrightarrow f(x)\sim f(y).
$$
If $(X,f)$ is   topologically mixing and $X/_\sim$ is metrizable, then $(X/_\sim,f^{\star})$ is   topologically mixing.
\end{theorem}
\begin{proof}
	See \cite[Theorem 3.22]{USS}
\end{proof}
\begin{theorem}\label{lijk}
	 The following holds. 
	 \begin{enumerate}
	 	\item The mapping $(\sigma_H^+)^{\star}$ is a non-invertible mapping on $\mathbb X_H^+/_{\sim_1}$ such that the dynamical system $(\mathbb X_H^+/_{\sim_1},(\sigma_H^+)^{\star})$ is topologically mixing.
	 	\item The mapping $\sigma_H^{\star}$ is a homeomorphism on $\mathbb X_H/_{\sim_2}$ such that the dynamical system $(\mathbb X_H/_{\sim_2},\sigma_H^{\star})$ is topologically mixing.
	 \end{enumerate} 
\end{theorem}
\begin{proof}
	By Theorem \ref{lejgaungatamlehre}, the dynamical systems $(\mathbb X_H^+,\sigma_H^+)$ and $(\mathbb X_H,\sigma_H)$ are topologically mixing. It follows from Theorem \ref{tupki} that also the dynamical systems $(\mathbb X_H^+/_{\sim_1},(\sigma_H^+)^{\star})$ and $(\mathbb X_H/_{\sim_2},\sigma_H^{\star})$ are topologically mixing.
\end{proof}

\begin{theorem}\label{lilija}
	There is a sequence $(A_n)$ in $\mathcal L(L)$ such that $\displaystyle \lim_{n\to \infty}\diam(A_n)=0$ and there is 
	\begin{enumerate}
		\item a non-invertible mapping $f:L\rightarrow L$ such that
		\begin{enumerate}
		\item  $(L,f)$ is topologically mixing and  
		\item for each positive integer $n$ and for each $x\in A_n$, $f(x)=x$.
	\end{enumerate}
	\item a homeomorphism $\varphi:L\rightarrow L$ such that 
	\begin{enumerate}
		\item  $(L,\varphi)$ is topologically mixing and  
		\item for each positive integer $n$ and for each $x\in A_n$, $\varphi(x)=x$.
	\end{enumerate}
	\end{enumerate} 
	\end{theorem}
\begin{proof}
	First, we prove that there are a sequence $(A_n)$ in $\mathcal L(L)$  and a non-invertible mapping $f:L\rightarrow L$ such that
		\begin{enumerate}
		\item $\displaystyle \lim_{n\to \infty}\diam(A_n)=0$,
		\item $(L,f)$ is topologically mixing, and  
		\item for each positive integer $n$ and for each $x\in A_n$, $f(x)=x$.
	\end{enumerate}
	For each positive integer $n$, let
	$$
	B_n=\{[(x_1,x_2,x_3,\ldots)]_1\in \mathbb X_H^+/_{\sim_1} \ | \  \mathbf x_1\in I_n, \textup{for all positive integers } k,\ell, \mathbf x(k)=\mathbf x(\ell)\}.
	$$
	Note that for each positive integer $n$,
	\begin{enumerate}
		\item $B_n\in \mathcal L(\mathbb X_H^+/_{\sim_1})$ and 
		\item for each $\mathbf x\in \mathbb X_H^+$, $(\sigma_H^+)^{\star}([\mathbf x]_1)=[\mathbf x]_1$.
	\end{enumerate} 
	Let $h:\mathbb X_H^+/_{\sim_1}\rightarrow L$ be a homeomorphism and let $f:L\rightarrow L$ be defined by 
	$$
	f(x)=h((\sigma_H^+)^{\star}(h^{-1}(x)))
	$$
	for each $x\in L$. Also, for each positive integer $n$, let $A_n=h(B_n)$. Then $f$ is a non-invertible map and $(A_n)$ is a sequence in $\mathcal L(L)$ such that 
	\begin{enumerate}
		\item $\displaystyle \lim_{n\to \infty}\diam(A_n)=0$,
		\item $(L,f)$ is topologically mixing, and  
		\item for each positive integer $n$ and for each $x\in A_n$, $f(x)=x$.
	\end{enumerate}
	The proof of the second part of the theorem is analogous to the proof of the first part. We leave the details to a reader.
\end{proof}

Theorem \ref{main1} is the main result of the paper.

\begin{theorem}\label{main1}
	There is an uncountable set $\Lambda$ and 
\begin{enumerate}
	\item a family 
	$$
	\mathcal F=\{(X_\lambda,f_\lambda) \ | \ \lambda \in \Lambda\}
	$$ 
	of dynamical systems  such that 
	\begin{enumerate}
		\item for each $\lambda\in \Lambda$, $X_\lambda$ is a Lelek-like fan,
		\item for each $\lambda\in \Lambda$, $\varphi_\lambda$ is a non-invertible mapping on $X_\lambda$ such that $(X_\lambda,f_\lambda)$ is topologically mixing, and
		\item for all $\lambda_1,\lambda_2\in \Lambda$,
		$$
		\lambda_1\neq \lambda_2 ~~~ \Longrightarrow ~~~   X_{\lambda_1} \textup{ is not homeomorphic to } X_{\lambda_2}.
		$$ 
	\end{enumerate}
	\item a family 
	$$
	\mathcal F=\{(X_\lambda,\varphi_\lambda) \ | \ \lambda \in \Lambda\}
	$$ 
	of dynamical systems  such that 
	\begin{enumerate}
		\item for each $\lambda\in \Lambda$, $X_\lambda$ is a Lelek-like fan,
		\item for each $\lambda\in \Lambda$, $\varphi_\lambda$ is a homeomorphism on $X_\lambda$ such that $(X_\lambda,\varphi_\lambda)$ is topologically mixing, and
		\item for all $\lambda_1,\lambda_2\in \Lambda$,
		$$
		\lambda_1\neq \lambda_2 ~~~ \Longrightarrow ~~~   X_{\lambda_1} \textup{ is not homeomorphic to } X_{\lambda_2}.
		$$ 
	\end{enumerate}
\end{enumerate}	
\end{theorem}
\begin{proof}
	To prove the first part of the theorem, let $f:L\rightarrow L$ be a non-invertible mapping (from the proof of Theorem \ref{lilija}) and let $(A_n)$ be a sequence in $\mathcal L(L)$ such that \begin{enumerate}
		\item $\displaystyle \lim_{n\to \infty}\diam(A_n)=0$,
		\item $(L,f)$ is topologically mixing, and  
		\item for each positive integer $n$ and for each $x\in A_n$, $f(x)=x$.
	\end{enumerate}
Without loss of generality, we assume that for each positive integer $n$, $A_n=A^n$. Let $\Lambda = \mathbf I$. Note that $\Lambda$ is an uncountable set. For each $\lambda\in \Lambda$, let $(X_{\lambda},f_{\lambda})$ be the quotient of the dynamical system $(L,f)$ that is obtained from the relation $\sim_{\lambda}$, and let 
	$$
	\mathcal F=\{(X_\lambda,f_\lambda) \ | \ \lambda \in \Lambda\}.
	$$ 
	Note that 
	\begin{enumerate}
		\item for each $\lambda\in \Lambda$, $X_\lambda$ is a Lelek-like fan,
		\item for each $\lambda\in \Lambda$, $f_\lambda$ is a non-invertible mapping on $X_\lambda$, and
		\item for all $\lambda_1,\lambda_2\in \Lambda$,
		$$
		\lambda_1\neq \lambda_2 ~~~ \Longrightarrow ~~~   X_{\lambda_1} \textup{ is not homeomorphic to } X_{\lambda_2}.
		$$ 
	\end{enumerate}
	The proof of the second part of the theorem is analogous to the proof of the first part. We leave the details to a reader.
	\end{proof}

	\section{Acknowledgement}
This work is supported in part by the Slovenian Research Agency (research projects J1-4632, BI-HR/23-24-011, BI-US/22-24-086 and BI-US/22-24-094, and research program P1-0285). 
	
%\section{Declarations}
%The following sections are not relevant to our manuscript. 
%\subsection{Competing interests}
%Not applicable.
%\subsection{Data Availability Statement}
%Not applicable.

\noindent I. Bani\v c\\
              (1) Faculty of Natural Sciences and Mathematics, University of Maribor, Koro\v{s}ka 160, SI-2000 Maribor,
   Slovenia; \\(2) Institute of Mathematics, Physics and Mechanics, Jadranska 19, SI-1000 Ljubljana, 
   Slovenia; \\(3) Andrej Maru\v si\v c Institute, University of Primorska, Muzejski trg 2, SI-6000 Koper,
   Slovenia\\
             {iztok.banic@um.si}           %  \\
     
				\-
				
		\noindent G.  Erceg\\
             Faculty of Science, University of Split, Rudera Bo\v skovi\' ca 33, Split,  Croatia\\
%             {i}     
{{goran.erceg@pmfst.hr}       }    %  \\

                 	\-
                 				
		\noindent I.  Jeli\' c\\
             Faculty of Science, University of Split, Rudera Bo\v skovi\' ca 33, Split,  Croatia\\
%             {i}     
{{ivajel@pmfst.hr}       }    %  \\

                 	\-
					
  \noindent J.  Kennedy\\
             Department of Mathematics,  Lamar University, 200 Lucas Building, P.O. Box 10047, Beaumont, Texas 77710 USA\\
%             {}     
{{kennedy9905@gmail.com}       }  

	\-
				

                 %  \\

%``text''
%%%%%%%%%%%%%%%%%%%%%%%%%%%%%%%%%%%%%%%%%%%%%%%%%%%%%%%%%%%%%%%%%%%%%%%%%%%%%%%%%
%%% I N T R O D U C T I O N S


\begin{thebibliography}{9}
%\bibitem{A} E. ~Akin, {General Topology of Dynamical Systems}, Volume 1, Graduate Studies in Mathematics Series, American Mathematical Society, Providence RI, 1993.
\bibitem{IJGI} I.~Bani\v{c}, G.~Erceg, I.~Jeli\' c, J.~Kennedy, V.~Nall,  Fans, phans and pans, 
https://doi.org/10.48550/arXiv.2504.07267.
\bibitem{BE} I.~Bani\v c, G.~Erceg,   J.~Kennedy, C.~Mouron, V.~Nall, Transitive mappings on the Cantor fan, Ergodic Theory and Dynamical Systems (2025) 1 -- 35,  https://doi.org/10.1017/etds.2025.6.
\bibitem{USS} I.~Bani\v c, G.~Erceg,   J.~Kennedy, C.~Mouron, V.~Nall, Chaos and mixing homeomorphisms on fans, J. Difference Equ. Appl. (2024) 1--31, https://doi.org/10.1080/10236198.2024.2384947.
\bibitem{short} I.~Bani\v c, G.~Erceg,   J.~Kennedy, An embedding of the Cantor fan into the Lelek fan, Rad HAZU 29 = 564 (2025) 221-229, https://doi.org/10.21857/m8vqrt3lk9.
\bibitem{BE2} I. Bani\v c, G. Erceg, J. Kennedy, Closed relations with non-zero entropy that generate no periodic points, Discrete Contin. Dyn. Syst., 42 (2022), 5137--5166.
\bibitem{banic1} I.~Bani\v c, G.~Erceg,  J.~Kennedy, The Lelek fan as the inverse limit of intervals with a single set-valued bonding function whose graph is an arc,  {Mediterr. J. Math.  20 (2023) 1--24}..
\bibitem{banic2} I.~Bani\v c, G.~Erceg,  J.~Kennedy, A transitive homeomorphism on the Lelek fan,  J. Difference Equ. Appl. 23 (2023) 393--418. %https://doi.org/10.1080/10236198.2023.2208242.
\bibitem{van} I.~Bani\v c, R.~Gril Rogina,  J.~Kennedy, V.~Nall, Sufficient conditions for non-zero entropy of closed relations, Ergodic Theory and Dynamical Systems (2024) 1--29, https://doi.org/10.1017/etds.2024.11.

\bibitem{borsuk1} K. Borsuk, \emph{A theorem on fixed points}, Bul. Acad. Bolon. Sci. Cl. 2 (1954) 17--20.

%\bibitem{banks} J. ~Banks, J. ~Brooks, G. ~Cairns, G. ~Davis and P. ~Stacey, On Devaney's Definition of Chaos, The American Mathematical Monthly 99 (1992) 332--334.
%\bibitem{devaney} R.~L.~Devaney, A first course in chaotic dynamical systems: theory and experiments. Massachusetts: Perseus Books, 1992.
\bibitem{oversteegen} W.~D.~Bula and L.~Overseegen, A Characterization of smooth Cantor Bouquets,  Proc. Amer.Math.Soc. 108 (1990) 529--534.
\bibitem{jjc1} J. J. Charatonik, {On ramification points on the classical sense}, Fund. Math. 51 (1962/1963) 229--252.
\bibitem{charatonik} W.~J.~Charatonik, The Lelek fan is unique, Houston J. Math. 15 (1989) 27--34.
\bibitem{Jcharatonik} J.~J. ~Charatonik,  On fans,  Dissertationes Math.  54 (1967).
\bibitem{Jcharatonik1} J.~J. ~Charatonik, W.~J.~Charatonik, Smoothness and the property of Kelley,  Comment. Math. Univ. Carolin. 41 (2000) 123--132.
\bibitem{chenli} L.~Chen, E.~Li, Shadowing property for inverse limit spaces, Proc. Amer. Math. Soc. 115 (1992) 573--580. 
\bibitem{Yorke} E.~Coven, I.~Kan, J.~Yorke, Pseudo-orbit shadowing in the family of tent maps,  {Trans. Amer. Math. Soc. 308 (1988) 227--241}
\bibitem{koch} R.~J.~Koch, Arcs in partially ordered spaces, Pacific J.  Math.  20 (1959) 723--728.
\bibitem{KO} D.~Kwietniak, P.~Oprocha, A note on the average shadowing property for expansive maps, Topology and its Applications 159 (2012) 19--27.
\bibitem{eberhart} C. ~Eberhart, A note on smooth fans, Colloq.  Math.  20 (1969) 89--90.
\bibitem{engelking1} R.~ Engelking, General topology, Heldermann, Berlin, 1989.
\bibitem{nagata} R. Engelking, \emph{Dimension theory}, North-Holland, Amsterdam, 1978.
\bibitem{KS} S.~Kolyada, L.~Snoha,  Topologically transitivity, \textit{Scholarpedia}  4 (2):5802 (2009).
%\bibitem{knudsen} C.~Knudsen, Chaos Without Nonperiodicity, The American Mathematical Monthly 101(1994) 563--565. 
\bibitem{lelek} A.~Lelek, On plane dendroids and their end-points in the classical sense, Fund. Math. 49 (1960/1961) 301--319.
\bibitem{mackowjak} T.~Ma\' ckowiak, The hereditary classes of mappings, Fund. Math. 97 (1977), 123--150.
%\bibitem{li} S. ~Li, Dynamical properties of the shift maps on the inverse limit spaces,  Ergod. Th.  Dynam. Sys.  12 (1992) 95--108.
\bibitem{nadler} S.~B.~Nadler, Continuum theory. An introduction, Marcel Dekker, Inc., New York, 1992.
\bibitem{piotr} P.~Oprocha, Lelek fan admits completely scrambled weakly topologically mixing homeomorphism,  Bulletin of the London Mathematical Society 57 (2025) 432--443.
\bibitem{piotrvan} P.~Oprocha, V.~Nall, Lelek fan admits a topologically mixing homeomorphism with any entropy, preprint, 2025.
\end{thebibliography}
\end{document}